\documentclass[12pt]{amsart}

\usepackage[T1]{fontenc}
\usepackage[utf8]{inputenc}
\usepackage{lmodern}
\usepackage{amsmath,amssymb,amsthm,mathtools}
\usepackage{mathrsfs}
\usepackage{hyperref}
\usepackage{enumitem}
\usepackage{microtype}
\usepackage{tikz}
\usepackage{float}
\usetikzlibrary{graphs}
\usetikzlibrary{graphs,quotes}
\usetikzlibrary{decorations.pathreplacing,decorations.markings}
\usetikzlibrary{decorations.pathmorphing}
\usepackage{longtable}

\usepackage[margin=1in]{geometry}

\numberwithin{equation}{section}
\hypersetup{colorlinks=true,linkcolor=blue,citecolor=blue,urlcolor=blue}

\newtheorem{theorem}{Theorem}[section]
\newtheorem{proposition}[theorem]{Proposition}
\newtheorem{conjecture}[theorem]{Conjecture}
\newtheorem{lemma}[theorem]{Lemma}
\newtheorem{corollary}[theorem]{Corollary}
\theoremstyle{definition}
\newtheorem{definition}[theorem]{Definition}

\theoremstyle{remark}

\title[]{A Proof of The Twin-Star Hypothesis in Dimension Five}
\author[]{ \tiny Patrick Nyadjo Fonga}

\address{Department of Mathematics and Statistics, Bowling Green State University,
Bowling Green, OH 43403, USA}
\email{nyadjop@bgsu.edu}

\date{}

\begin{document}

\begin{abstract}
It was recently shown that the question of whether the determinant-like map
\(\det^{S^2}\) is unique up to multiplication by a scalar
is equivalent to a combinatorial conjecture regarding
cycle-free \(d\)-partitions of the complete graph \(K_{2d}\), known as the twin-star hypothesis. This hypothesis asserts that every cycle-free \(d\)-partition of \(K_{2d}\) is weakly equivalent to one that contains a component isomorphic to the twin-star graph \(TS_d\). This hypothesis was previously established for \(d\le4\).
In the present paper, we prove it for \(d=5\). Consequently, the associated determinant-like map \(\det^{S^2}\) is unique up to a scalar in dimension \(5\).
\end{abstract}

\maketitle

\section{Introduction}

The decomposition of complete graphs into spanning trees is a
natural problem in graph theory, and has been investigated from
several perspectives. In particular, decompositions of complete
graphs of even order into isomorphic spanning trees have been
studied through graph factorizations and labeling methods
\cite{Froncek,KubesaSpanningTrees}.
Further work has examined which tree structures can occur in
such factorizations, including trees with a prescribed number
of non-leaf vertices and particular families of trees
\cite{FroncekKovarKubesa,KovarKubesaMeszka}.
Although the structure of an individual spanning tree is well
understood, requiring the edges of a complete graph to be
distributed among several spanning trees, with each edge
belonging to exactly one tree, imposes strong restrictions on
their possible structures. Related questions have been
studied under additional conditions, such as edge-coloring
constraints
\cite{BrualdiHollingsworth, GlockKuhnMontgomeryOsthus, MazzuoccoloRinaldi}
and geometric restrictions
\cite{BoseHurtadoRiveraWood}.
More general questions regarding partitions of complete
graphs have also been considered in \cite{CioabaCameron}.

Let $d\ge2$ be an integer, and let $K_{2d}$ denote the
complete graph on $2d$ vertices. A \emph{cycle-free
$d$-partition} of $K_{2d}$ is an ordered $d$-tuple
\[
\mathcal P=(\Gamma_1,\ldots,\Gamma_d)
\]
of cycle-free spanning subgraphs whose edge sets are
pairwise disjoint and together contain every edge of
$K_{2d}$. Since $K_{2d}$ has $d(2d-1)$ edges and
each cycle-free spanning subgraph has at most $2d-1$
edges, every component $\Gamma_i$ must have exactly
$2d-1$ edges and is therefore a spanning tree.

Our interest is in understanding how these partitions
can be transformed into one another while preserving
the cycle-free property. A particularly useful family
of transformations was introduced by Lippold, Staic,
and Stancu in \cite{edge}. Given a cycle-free
$d$-partition $\mathcal P$ and three distinct vertices
$x,y,z$ of $K_{2d}$, they proved that there exists
a unique cycle-free $d$-partition
$\mathcal P^{(x,y,z)}$ which agrees with $\mathcal P$
on every edge outside the triangle with vertices
$x,y,z$, but differs from $\mathcal P$ on at least
two of the three edges of this triangle. Moreover,
the transformation
\[
\mathcal P\longmapsto\mathcal P^{(x,y,z)}
\]
is an involution, called the \emph{triangle involution
$(x,y,z)$}, or simply the \emph{involution $(x,y,z)$}.

Triangle involutions provide a way to modify a
cycle-free partition locally without destroying
the spanning-tree structure of its components. By applying these
transformations, we obtain a natural
equivalence relation on the set of cycle-free
$d$-partitions of $K_{2d}$.

Two cycle-free $d$-partitions are said to be
\emph{involution equivalent} if one can be obtained
from the other by a finite sequence of triangle
involutions. This leads to the following conjecture, formulated
in \cite{edge}.

\begin{conjecture}\label{conj:Cd}
Any two cycle-free $d$-partitions of $K_{2d}$
are involution equivalent.
\end{conjecture}

Conjecture \ref{conj:Cd} was established for $d=3$ by
Lippold, Staic, and Stancu \cite{edge}. A further reduction of this conjecture is obtained in \cite{twin}. The Authors showed that the problem can be approached by focusing on a particular spanning tree among the
components of a partition. To formulate it, we need to recall the following definition from \cite[Definition 3.9]{twin}.

\begin{definition}\label{def:weak-equivalence}
Two cycle-free $d$-partitions of $K_{2d}$ are said
to be \emph{weakly equivalent} if, after possibly
relabeling the vertices of $K_{2d}$ and permuting
the components of one partition, the resulting
partitions are involution equivalent.
\end{definition}

\begin{figure}[ht]
\centering
\begin{tikzpicture}[scale=0.76,
  vertex/.style={circle,draw=black,minimum size=8mm,inner sep=0pt},
  edge/.style={line width=0.45mm}]
\node[vertex] (one) at (-1.35,0) {$1$};
\node[vertex] (two) at (1.35,0) {$2$};
\draw[edge] (one)--(two);
\node[vertex] (l1) at (-1.80,2.40) {$3$};
\node[vertex] (l2) at (-3.62,2.35) {$5$};
\node[vertex] (l3) at (-5.12,1.00) {$7$};
\node[vertex] (l4) at (-5.12,-0.98) {\scriptsize $2d-5$};
\node[vertex] (l5) at (-3.62,-2.35) {\scriptsize $2d-3$};
\node[vertex] (l6) at (-1.80,-2.40) {\scriptsize $2d-1$};
\node[vertex] (r1) at (1.80,2.40) {$4$};
\node[vertex] (r2) at (3.62,2.35) {$6$};
\node[vertex] (r3) at (5.12,1.00) {$8$};
\node[vertex] (r4) at (5.12,-0.98) {\scriptsize $2d-4$};
\node[vertex] (r5) at (3.62,-2.35) {\scriptsize $2d-2$};
\node[vertex] (r6) at (1.80,-2.40) {\scriptsize $2d$};
\foreach \u in {l1,l2,l3,l4,l5,l6} \draw[edge] (one)--(\u);
\foreach \u in {r1,r2,r3,r4,r5,r6} \draw[edge] (two)--(\u);
\node at (-5.12,0) {$\vdots$};
\node at (5.12,0) {$\vdots$};
\end{tikzpicture}
\caption{The twin-star graph $TS_d$, with $2d$ vertices.}
\label{fig:tsd}
\end{figure}
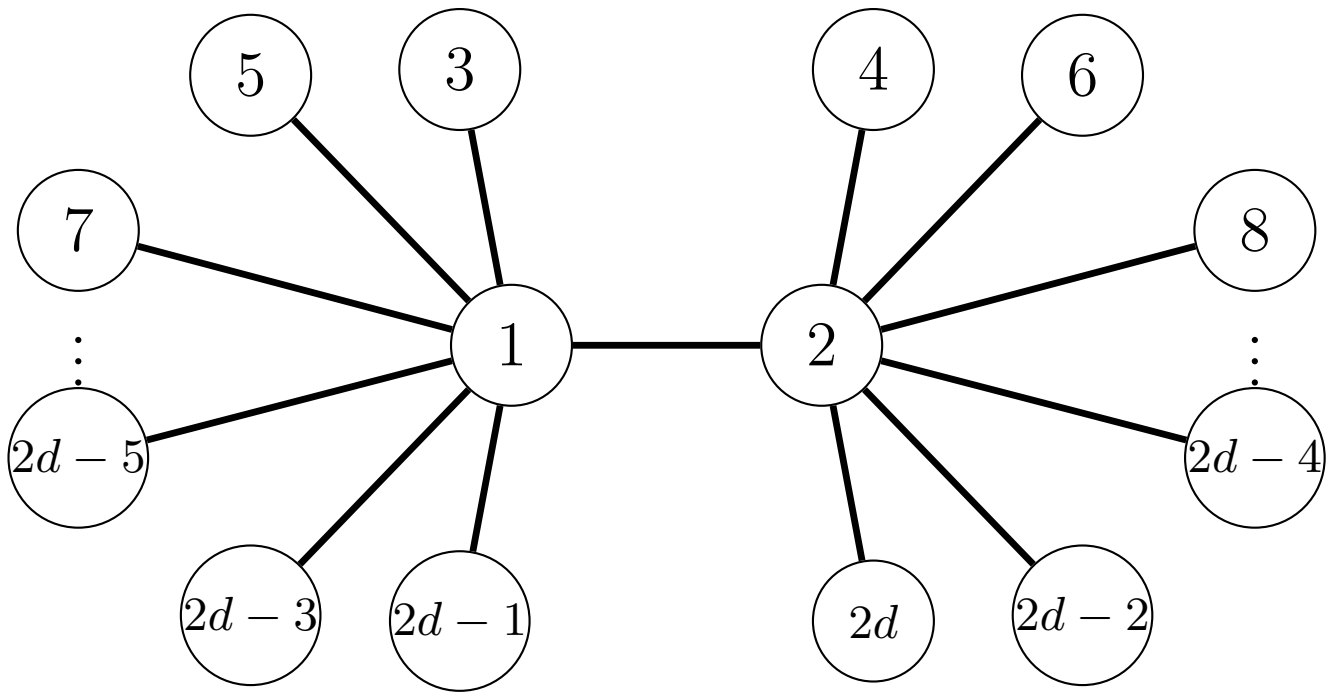

The particular spanning tree targeted in this reduction
is the \emph{twin-star graph} $TS_d$. It consists
of two adjacent vertices of degree $d$, each
joined to $d-1$ leaves. Thus $TS_d$ is a tree
on $2d$ vertices with exactly $2d-2$ leaves;
see Figure~\ref{fig:tsd}. 
With this terminology, the reduction of
Conjecture~\ref{conj:Cd}, known as the
\emph{twin-star hypothesis}, can be stated as follows.

\begin{conjecture}
\label{conj:TSd}
Every cycle-free $d$-partition of $K_{2d}$
is weakly equivalent to one containing a component
isomorphic to $TS_d$.
\end{conjecture}

The relationship between Conjectures~\ref{conj:Cd}
and~\ref{conj:TSd} was established in \cite{twin}.
If Conjecture~\ref{conj:Cd} holds, then
Conjecture~\ref{conj:TSd} follows. Conversely,
if the twin-star hypothesis holds for every
$1\leq i\leq d$, then Conjecture~\ref{conj:Cd}
also holds. Using this reduction, the authors
of \cite{twin} established the case $d=4$,
extending the previously known cases $d=1$, $d=2$
and $d=3$. The advantage of this formulation is that it reduces the problem of comparing arbitrary
cycle-free partitions to that of obtaining a
specific tree structure in one of their components.

These combinatorial questions also have an
important algebraic motivation. They arise from a uniqueness problem
concerning a determinant-like map associated with
the generalized exterior algebra $\Lambda_{V_d}^{S^2}$
introduced by Staic in \cite{sta2}.
Let $k$ be a field and let $V_d$ be a
$d$-dimensional vector space over $k$.
Motivated by a generalization of the classical
determinant, Staic conjectured in \cite{sta2}
the existence and uniqueness, up to multiplication
by a scalar, of a nonzero linear map
\[
\mathrm{det}^{S^2}:V_d^{\otimes d(2d-1)}
\longrightarrow k
\]
satisfying the following vanishing condition:
\[
\mathrm{det}^{S^2}\left(
\bigotimes_{1\leq i<j\leq 2d}v_{i,j}
\right)=0
\]
whenever there exist three indices
\(1\leq x<y<z\leq2d\) such that
\[
v_{x,y}=v_{x,z}=v_{y,z}.
\]
The connection
between this algebraic construction and the
combinatorial structure of $K_{2d}$ was developed
in \cite{edge}, where the uniqueness problem
for $\det^{S^2}$ was shown to be equivalent to
Conjecture \ref{conj:Cd}.

The existence of $\det^{S^2}$ for arbitrary
$d$ was established in \cite{dets2}, while its
uniqueness in general remains open. The map
also has applications in geometry and
combinatorics \cite{dets2,sv}, as well as in
equilibrium problems in classical mechanics
\cite{req}. Consequently, progress on the
twin-star hypothesis also provides information
about the uniqueness of this determinant-like map.

The purpose of the present paper is to establish
the twin-star hypothesis in dimension five. When $d=5$, a cycle-free $5$-partition $\mathcal P$ of
$K_{10}$ consists of five spanning trees, each having
$9$ edges. The twin-star $TS_5$ is the tree on $10$
vertices with two adjacent vertices of degree $5$
and eight leaves. We prove that every cycle-free
$5$-partition of $K_{10}$ is weakly equivalent to one
containing a component isomorphic to $TS_5$, thereby
establishing Conjecture~\ref{conj:TSd} for $d=5$.
Since the twin-star hypothesis is already known for
$1\leq d\leq4$, our result also establishes
Conjecture~\ref{conj:Cd} for $d=5$.
Consequently, the corresponding map
\[
\mathrm{det}^{S^2}:V_5^{\otimes45}\longrightarrow k
\]
is unique up to multiplication by a scalar.

The paper is organized as follows. In Section~\ref{sec:prelim},
we introduce the notation and establish some preliminary results
on cycle-free partitions. In Section~\ref{sec:extremal}, we give more properties of cycle-free $d$-partitions
and introduce the notion of an extremal component, which is a component chosen
to have the maximum number of leaves among all the component of all partition weakly equivalent to a given partition, and among all such
components, the minimum number of branch vertices.
Using weak equivalence, we establish the structural
properties of this component and prove that it has
at least $d+1$ leaves; see
Theorem~\ref{thm:general-lower-bound}.
In Section~\ref{sec:d5-prelim}, we specialize to $d=5$
and determine the possible branch-degree configurations
of an extremal component with six or seven leaves.
Sections~\ref{sec:d5-six} and~\ref{sec:d5-seven}
exclude these two cases, respectively, showing that
an extremal component must have eight leaves and
hence be isomorphic to $TS_5$. This establishes
our main result in Theorem~\ref{thm:main}.

\section{Preliminaries}\label{sec:prelim}

In this section, we introduce the notation used throughout
the paper and establish some preliminary results
on cycle-free partitions that we use in the subsequent sections.

Let $H$ be a graph. We denote its vertex and edge sets
by $V(H)$ and $E(H)$, respectively. An edge joining
two vertices $u$ and $v$ is denoted by $(u,v)$.
We write $d_H(v)$ for the degree of a vertex $v$
in $H$ and $N_H(v)$ for the set of its neighbors. We also use the notation
\[
u_0-u_1-\cdots-u_m
\]
to denote a path with consecutive vertices
$u_0,u_1,\ldots,u_m$.

We refer to vertices of degree at least $3$ as
\emph{branch vertices} and we use the notations
\[
\begin{aligned}
\ell(H)&=\text{number of leaves of }H,\\
p(H)&=\text{number of vertices of degree }2 \ \text{of} \ H,\\
b(H)&=\text{number of branch vertices of }H.
\end{aligned}
\]

Since every vertex of a tree with at least two vertices
is either a leaf, a degree-$2$ vertex, or a branch
vertex, we have
\begin{equation}\label{eq:count-lpb}
\ell(H)+p(H)+b(H)=|V(H)|.
\end{equation}

We recall the following identity relating
the number of leaves of a tree to the degrees of
its branch vertices.

\begin{lemma}\label{lem:leaf-identity}
For every tree \(H\) with at least two vertices,
\[
\ell(H)=2+\sum_{d_H(v)\ge3}\bigl(d_H(v)-2\bigr).
\]
\end{lemma}

\begin{proof}
Since \(H\) is a tree, we have that
\[
\sum_{v\in V(H)}d_H(v)=2(|V(H)|-1),
\]
and consequently
\[
\sum_{v\in V(H)}\bigl(d_H(v)-2\bigr)=-2.
\]
Splitting the sum according to vertex degree, we obtain
\[
-\ell(H)+\sum_{d_H(v)\ge3}\bigl(d_H(v)-2\bigr)=-2,
\]
which gives the desired identity.
\end{proof}

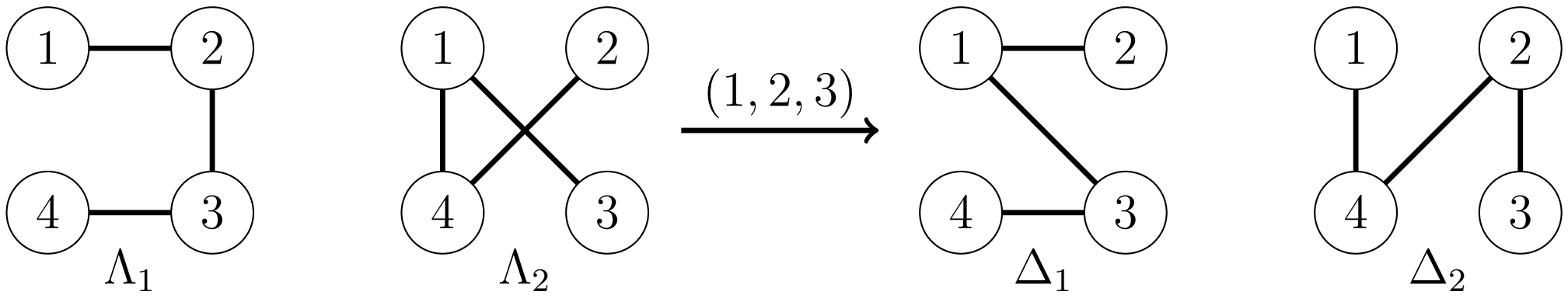
\begin{figure}[ht]
\centering
\begin{tikzpicture}[
    scale=0.7,
    vertex/.style={
        circle,draw=black,fill=white,
        minimum size=7mm,inner sep=0pt
    },
    edge/.style={line width=0.45mm}
]

\begin{scope}[shift={(0,0)}]
\coordinate (a) at (0,2);
\coordinate (b) at (2,2);
\coordinate (c) at (2,0);
\coordinate (d) at (0,0);

\draw[edge] (a)--(b)--(c)--(d);

\foreach \p/\n in {a/1,b/2,c/3,d/4}
    \node[vertex] at (\p) {$\n$};

\node at (1,-0.7) {$\Lambda_1$};
\end{scope}

\begin{scope}[shift={(4.8,0)}]
\coordinate (a) at (0,2);
\coordinate (b) at (2,2);
\coordinate (c) at (2,0);
\coordinate (d) at (0,0);

\draw[edge] (a)--(c);
\draw[edge] (a)--(d);
\draw[edge] (b)--(d);

\foreach \p/\n in {a/1,b/2,c/3,d/4}
    \node[vertex] at (\p) {$\n$};

\node at (1,-0.7) {$\Lambda_2$};
\end{scope}

\draw[->,line width=0.45mm] (7.7,1)--(10.1,1);
\node at (8.9,1.45) {$(1,2,3)$};

\begin{scope}[shift={(11.1,0)}]
\coordinate (a) at (0,2);
\coordinate (b) at (2,2);
\coordinate (c) at (2,0);
\coordinate (d) at (0,0);

\draw[edge] (a)--(b);
\draw[edge] (a)--(c);
\draw[edge] (c)--(d);

\foreach \p/\n in {a/1,b/2,c/3,d/4}
    \node[vertex] at (\p) {$\n$};

\node at (1,-0.7) {$\Delta_1$};
\end{scope}

\begin{scope}[shift={(15.9,0)}]
\coordinate (a) at (0,2);
\coordinate (b) at (2,2);
\coordinate (c) at (2,0);
\coordinate (d) at (0,0);

\draw[edge] (b)--(c);
\draw[edge] (a)--(d);
\draw[edge] (b)--(d);

\foreach \p/\n in {a/1,b/2,c/3,d/4}
    \node[vertex] at (\p) {$\n$};

\node at (1,-0.7) {$\Delta_2$};
\end{scope}

\end{tikzpicture}
\caption{An illustration of the triangle involution $(1,2,3)$
on the cycle-free $2$-partition $\mathcal T=(\Lambda_1,\Lambda_2)$ of $K_4$.}
\label{fig:triangle-example}
\end{figure}

Before studying the effect of triangle involutions on cycle-free
$d$-partitions of $K_{2d}$,
we first illustrate how such an involution changes the
distribution of edges among the components of a cycle-free
partition in general. Consider as an example the cycle-free $2$-partition
\(
\mathcal T=(\Lambda_1,\Lambda_2)
\)
of $K_4$ defined by
\[
\begin{aligned}
E(\Lambda_1)=\{(1,2),(2,3),(3,4)\}, \quad \text{and} \quad
E(\Lambda_2)=\{(1,3),(1,4),(2,4)\},
\end{aligned}
\]
as shown in Figure~\ref{fig:triangle-example}.

We apply the involution $(1,2,3)$. Since the involution
changes only the three edges of the triangle with vertices
$1,2,3$, the only possible nontrivial changes involve the
edge $(1,3)$ and one of the edges $(1,2)$ or $(2,3)$.
The edge $(1,3)$ cannot be exchanged with $(1,2)$, because
this would give the second component the edge set
\[
\{(1,2),(1,4),(2,4)\},
\]
which contains a cycle. Therefore this possibility is excluded.
Since the involution is nontrivial, it follows that it must
exchange the edges $(2,3)$ and $(1,3)$ between the two components.

Thus, writing
\(
\mathcal T^{(1,2,3)}=(\Delta_1,\Delta_2),
\)
we obtain
\[
\begin{aligned}
E(\Delta_1)=\{(1,2),(1,3),(3,4)\},\quad \text{and} \quad
E(\Delta_2)=\{(2,3),(1,4),(2,4)\}.
\end{aligned}
\]
Both resulting components are again spanning trees.

Let $d\ge2$ and let
$\mathcal P=(\Gamma_1,\ldots,\Gamma_d)$ be a cycle-free
$d$-partition of $K_{2d}$. The fact that the components of a cycle-free
$d$-partition together contain every edge of
$K_{2d}$ imposes a restriction on the
degree of each vertex.

\begin{lemma}\label{lem:degree-bound}
For every vertex $v\in V(K_{2d})$ and every
$1\le i\le d$, we have
\[
1\le d_{\Gamma_i}(v)\le d.
\]
\end{lemma}

\begin{proof}
The $2d-1$ edges incident with $v$ are distributed
among the $d$ spanning trees
$\Gamma_1,\ldots,\Gamma_d$, so
\[
\sum_{i=1}^d d_{\Gamma_i}(v)=2d-1.
\]
Since each $\Gamma_i$ is connected, every summand \(d_{\Gamma_i}(v)\)
is at least $1$. Hence, for any fixed
$j\in\{1,\ldots,d\}$,
\[
d_{\Gamma_j}(v)
\le(2d-1)-(d-1)=d.
\]
\end{proof}

The following lemma describes how
an involution can change the attachment of a leaf in a component of the cycle-free
$d$-partition $\mathcal P$.

\begin{lemma}\label{lem:leaf-switch}
Let $v$ be a leaf of $\Gamma_d$, with
$(a,v)\in E(\Gamma_d)$, and suppose that
$(a,b)$ and $(b,v)$ do not belong to
$E(\Gamma_d)$. Then the involution
$(a,b,v)$ has one of the following two effects:
\begin{enumerate}[label=\textup{(\roman*)}]
\item It replaces $(a,v)$ by $(b,v)$ in
$\Gamma_d$.
\item It leaves $\Gamma_d$ unchanged, in which case
$(a,b)$ and $(b,v)$ belong to two different
components and are exchanged by the involution.
\end{enumerate}
In particular, if $(a,b)$ and $(b,v)$ belong
to the same component, then the involution
replaces $(a,v)$ by $(b,v)$ in $\Gamma_d$.
\end{lemma}

\begin{proof}
Among the three edges of the triangle with
vertices $a,b,v$, the component $\Gamma_d$
contains only $(a,v)$. If the involution changes
$\Gamma_d$, this edge must be replaced by
either $(a,b)$ or $(b,v)$. The first possibility
is impossible, since deleting the leaf edge
$(a,v)$ isolates $v$, while adding $(a,b)$
does not reconnect it. Thus the only possible
change in $\Gamma_d$ is the replacement
of $(a,v)$ by $(b,v)$.

If $\Gamma_d$ remains unchanged, the other
two edges must be exchanged between two distinct
components, because the involution is nontrivial.
This also proves the final assertion.
\end{proof}

The possibility (i) in Lemma~\ref{lem:leaf-switch}
allows us to move a leaf from one vertex to another
while preserving the spanning-tree structure.
We introduce the following terminology to
describe this type of operation.

\begin{definition}\label{def:leaf-transfer}
Let $H$ be a component of a cycle-free
$d$-partition, and let $v$ be a leaf of $H$
adjacent to a vertex $a$. For a vertex
$b\notin\{a,v\}$, if the involution $(a,b,v)$
replaces the edge $(a,v)$ by $(b,v)$ in $H$,
we say that it \emph{transfers the leaf $v$
from $a$ to $b$}, or that $b$ \emph{receives
the leaf $v$}. In this case, $b$ is called
a \emph{destination} of $v$.
\end{definition}

The following lemma provides a sufficient condition for the existence
of such a destination for a leaf in \(\Gamma_d\). The particular case \(d=4\) was established
in \cite{twin}.

\begin{lemma}\label{lem:mobility}
Suppose that \(v\) is a leaf of \(\Gamma_d\), with \((a,v)\in E(\Gamma_d)\), and let
\[
B\subseteq V(K_{2d})\setminus\{a,v\}
\]
satisfy \((a,b)\notin E(\Gamma_d)\) for every \(b\in B\). If
\[
|B|\ge d-d_{\Gamma_d}(a)+1,
\]
then there exists \(b\in B\) that receives the leaf \(v\).
\end{lemma}

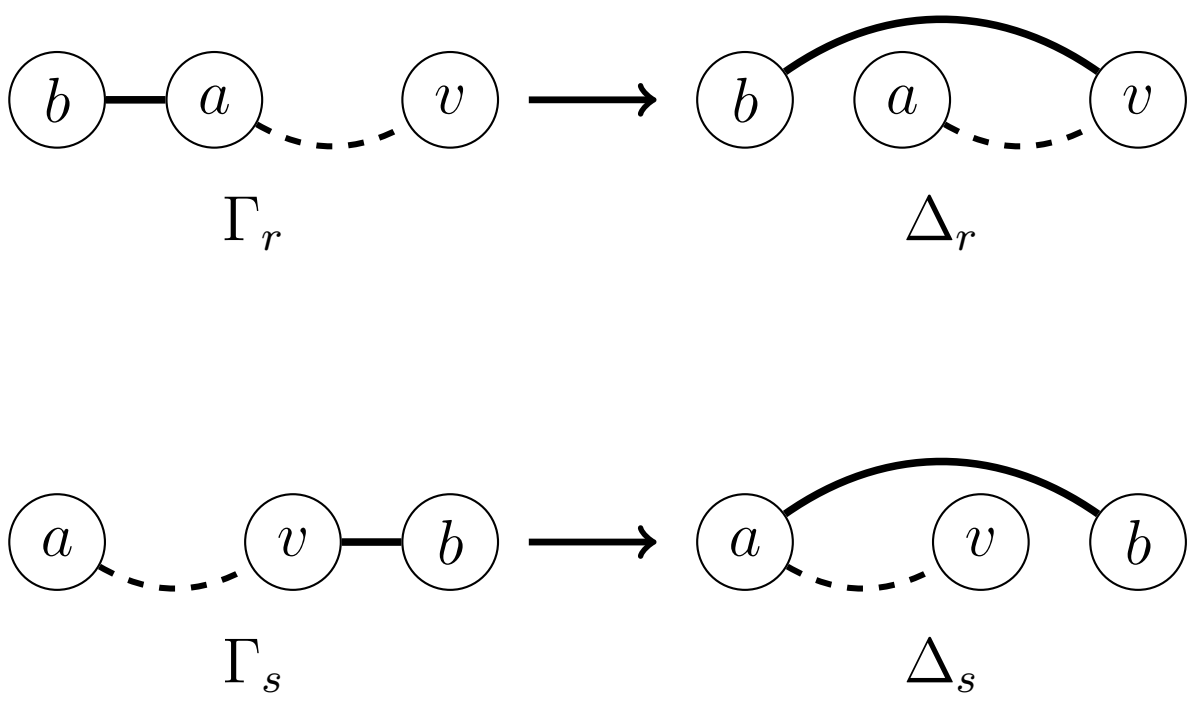
\begin{figure}[ht]
\centering
\begin{tikzpicture}[
    scale=0.65,
    vertex/.style={
        circle,draw=black,
        minimum size=18pt,
        inner sep=0.3pt
    }
]

\node[vertex] (b1) at (0,3) {$b$};
\node[vertex] (a1) at (1.6,3) {$a$};
\node[vertex] (v1) at (4,3) {$v$};

\draw[line width=0.5mm] (b1)--(a1);
\draw[dashed,line width=0.4mm]
    (a1) to[bend right=30] (v1);

\node at (2,1.75) {$\Gamma_r$};

\draw[->,line width=0.45mm] (4.8,3)--(6.1,3);

\node[vertex] (b2) at (7,3) {$b$};
\node[vertex] (a2) at (8.6,3) {$a$};
\node[vertex] (v2) at (11,3) {$v$};

\draw[line width=0.5mm]
    (b2) to[bend left=35] (v2);
\draw[dashed,line width=0.4mm]
    (a2) to[bend right=30] (v2);

\node at (9,1.75) {$\Delta_r$};

\node[vertex] (a3) at (0,-1.5) {$a$};
\node[vertex] (v3) at (2.4,-1.5) {$v$};
\node[vertex] (b3) at (4,-1.5) {$b$};

\draw[dashed,line width=0.4mm]
    (a3) to[bend right=30] (v3);
\draw[line width=0.5mm] (v3)--(b3);

\node at (2,-2.75) {$\Gamma_s$};

\draw[->,line width=0.45mm]
    (4.8,-1.5)--(6.1,-1.5);

\node[vertex] (a4) at (7,-1.5) {$a$};
\node[vertex] (v4) at (9.4,-1.5) {$v$};
\node[vertex] (b4) at (11,-1.5) {$b$};

\draw[dashed,line width=0.4mm]
    (a4) to[bend right=30] (v4);
\draw[line width=0.5mm]
    (a4) to[bend left=35] (b4);

\node at (9,-2.75) {$\Delta_s$};

\end{tikzpicture}
\caption{The exchange of $(a,b)$ and $(b,v)$
between $\Gamma_r$ and $\Gamma_s$.}
\label{fig:path-exchange}
\end{figure}

\begin{proof}
Suppose, for a contradiction, that no vertex \(b\in B\) receives the leaf \(v\). We fix \(b\in B\), and write
\[
(a,b)\in E(\Gamma_r),\quad \text{and} \quad
(b,v)\in E(\Gamma_s),
\]
for some \(r,s\in\{1,\ldots,d-1\}\). By Lemma~\ref{lem:leaf-switch}, we have \(r\ne s\), and we apply the involution \((a,b,v)\) which leaves \(\Gamma_d\) unchanged while exchanging \((a,b)\) and \((b,v)\).

Let \(\Delta_r\) denote the new \(r\)-th component. Then
\[
E(\Delta_r)=
\bigl(E(\Gamma_r)\setminus\{(a,b)\}\bigr)
\cup\{(b,v)\}.
\]
Since both \(\Gamma_r\) and \(\Delta_r\) are trees, the edge \((a,b)\) lies on the unique path from \(b\) to \(v\) in \(\Gamma_r\). Since \((a,b)\) is incident with \(b\), this path has the form
\[
b-a-\cdots-v.
\]
Consequently, the unique path from \(a\) to \(v\) in \(\Gamma_r\) does not start with \((a,b)\). This situation is illustrated in Figure~\ref{fig:path-exchange}.

Similarly, if \(\Delta_s\) denotes the new \(s\)-th component, then
\[
E(\Delta_s)=
\bigl(E(\Gamma_s)\setminus\{(b,v)\}\bigr)
\cup\{(a,b)\}.
\]
Since \(\Gamma_s\) and \(\Delta_s\) are trees, the edge \((b,v)\) belongs to the unique path from \(a\) to \(b\) in \(\Gamma_s\). In particular, the path from \(a\) to \(v\) in \(\Gamma_s\) does not start with \((a,b)\). The same is true for every other component \(\Gamma_i\), since none contains \((a,b)\).

Therefore, for every \(b\in B\), none of the \(d-1\) paths from \(a\) to \(v\) in \(\Gamma_1,\ldots,\Gamma_{d-1}\) can start with \((a,b)\).

There are \(2d-1\) edges incident with \(a\) and out of these, \(d_{\Gamma_d}(a)\) belong to \(\Gamma_d\), while the \(|B|\) edges \((a,b)\) with \(b\in B\), are excluded as possible first edges of the paths above. Hence these \(d-1\) paths have at most
\[
(2d-1)-d_{\Gamma_d}(a)-|B|
\]
possible first edges. Since the components are edge-disjoint, the paths must start with \(d-1\) distinct edges. Therefore
\[
(2d-1)-d_{\Gamma_d}(a)-|B|\ge d-1,
\]
which gives
\[
|B|\le d-d_{\Gamma_d}(a),
\]
contrary to the hypothesis.
\end{proof}

A consequence of Lemma \ref{lem:mobility} is the existence of a lower bound on the number of possible destinations of any leaf.

\begin{corollary}\label{cor:dminusone}
Let \(v\) be a leaf of \(\Gamma_d\), with \((a,v)\in E(\Gamma_d)\). Among the vertices \(b\notin\{a,v\}\) satisfying \((a,b)\notin E(\Gamma_d)\), at least \(d-1\) receive the leaf \(v\).
\end{corollary}

\begin{proof}

Put $m=d_{\Gamma_d}(a)$. Among the $2d-2$ vertices
different from $a$ and $v$, exactly $m-1$ are
adjacent to $a$ in $\Gamma_d$. Hence there are
\[
(2d-2)-(m-1)=2d-m-1
\]
vertices $b\notin\{a,v\}$ satisfying
$(a,b)\notin E(\Gamma_d)$.

We now determine how many of these vertices fail
to receive the leaf $v$. By Lemma~\ref{lem:mobility},
every set of at least $d-m+1$ such vertices contains
a destination of $v$. Therefore, the number of such
vertices that fail to receive $v$ is at most $d-m$.
Since there are $2d-m-1$ eligible vertices, at least
\[
(2d-m-1)-(d-m)=d-1
\]
of them receive the leaf $v$.
\end{proof}

We conclude this section by recalling the following
characterization of the twin-star graph $TS_d$,
established in \cite{twin}. We include a short proof
for the convenience of the reader.

\begin{lemma}\label{lem:twinstar-characterization}
Let $H$ be a component of a cycle-free
$d$-partition of $K_{2d}$. Then
\[
\ell(H)\le2d-2.
\]
Moreover, $\ell(H)=2d-2$ if and only if, $H$ is isomorphic to $TS_d$.
\end{lemma}

\begin{proof}
By Lemma~\ref{lem:degree-bound}, every vertex of $H$
has degree at most $d$. If $\ell(H)=2d-1$, then
$H$ is a star whose central vertex has degree
$2d-1>d$, which is a contradiction. Hence
\[
\ell(H)\le 2d-2.
\]

Suppose now that $\ell(H)=2d-2$. Then $H$
has exactly two non-leaf vertices, say $a$ and $b$.
Since $H$ is connected, these vertices must be
adjacent, and every other vertex is a leaf
adjacent to either $a$ or $b$. Consequently,
\[
d_H(a)+d_H(b)=2d.
\]
Since both degrees are at most $d$, we obtain
\[
d_H(a)=d_H(b)=d.
\]
Thus each of $a$ and $b$ is adjacent to
$d-1$ leaves, so $H$ is isomorphic to $TS_d$.

Conversely, the twin-star $TS_d$ has exactly
$2d-2$ leaves, which completes the proof.
\end{proof}

This characterization shows that the twin-star
hypothesis can be established by finding,
within each weak-equivalence class of
cycle-free $d$-partitions, a component with
$2d-2$ leaves.

\section{An extremal component}\label{sec:extremal}

In this section, we establish a general lower bound for the number of
leaves of a suitable component of a cycle-free $d$-partition. Our approach is to choose a component with as many leaves as possible and, among such components, with as few branch vertices as possible. We then use weak equivalence to rule out certain configurations, thereby restricting the structure of this component.

Throughout this section, let $d\ge3$ and let $\mathcal P$ be a cycle-free
$d$-partition of $K_{2d}$. We first give the following definition.

\begin{definition}\label{def:extremal-component}
A spanning tree $H$ is called an \emph{extremal component} of
$\mathcal P$ if it occurs as a component of a cycle-free $d$-partition
weakly equivalent to $\mathcal P$ and satisfies the following conditions:
\[
\begin{array}{ll}
\textup{(E1)} & \ell(H)\text{ is maximal among all components of all partitions weakly equivalent to $\mathcal P$},\\[2mm]
\textup{(E2)} & b(H)\text{ is minimal among the components satisfying \textup{(E1)}}.
\end{array}
\]
\end{definition}

Since there are only finitely many cycle-free $d$-partitions of
$K_{2d}$, an extremal component exists. Let \(T\) be one of such component and fix a weakly equivalent partition
$(\Gamma_1,\ldots,\Gamma_d)$ containing $T$. After permuting the
components if necessary, we may assume that $T=\Gamma_d$.

We now establish some properties of the component \(T\). The first shows that a degree-\(2\) vertex of \(T\) cannot lie between two non-leaf vertices.

\begin{lemma}\label{lem:deg2-adj-leaf}
Every degree-\(2\) vertex of the extremal component \(T\) is adjacent to a leaf.
\end{lemma}

\begin{proof}
Let \(v\in V(T)\) satisfy \(d_T(v)=2\), and write \(N_T(v)=\{x,y\}\). Suppose that neither \(x\) nor \(y\) is a leaf. Since \(T\) is a tree, \((x,y)\notin E(T)\). We apply the triangle involution \((v,x,y)\). On this triangle, the component \(T=\Gamma_d\) contains the edges \((v,x)\) and \((v,y)\), so after the involution one of these edges is replaced by \((x,y)\). Consequently, \(v\) becomes a leaf, while neither \(x\) nor \(y\) becomes a leaf. The new \(d\)-th component therefore has one more leaf than \(T\), contradicting \textup{(E1)}.
\end{proof}

An immediate consequence of Lemma \ref{lem:deg2-adj-leaf} is that degree-\(2\) vertices of \(T\) cannot be adjacent.

\begin{corollary}\label{cor:no-adjacent-deg2}
Two degree-\(2\) vertices of the extremal component \(T\) cannot be adjacent.
\end{corollary}

\begin{proof}
Suppose that \(u,v\in V(T)\) are adjacent and satisfy \(d_T(u)=d_T(v)=2\). We write
\[
N_T(u)=\{a,v\},\quad \text{and} \quad N_T(v)=\{u,b\}.
\]
By Lemma~\ref{lem:deg2-adj-leaf}, the vertices \(a\) and \(b\) are leaves. Hence \(a-u-v-b\) is a path in \(T\). Since \(a,b\) are leaves and \(u,v\) have degree \(2\), none of these four vertices is adjacent to any vertex outside this path. Thus the path is a connected component of \(T\), which is impossible because \(T\) is connected and has \(2d\ge6\) vertices.
\end{proof}

By Lemma~\ref{lem:deg2-adj-leaf} and
Corollary~\ref{cor:no-adjacent-deg2}, every degree-$2$ vertex of $T$
is adjacent to exactly one leaf and one branch vertex. In
particular, no degree-$2$ vertex of $T$ lies between two branch
vertices. This leads
to the following terminology.

\begin{definition}\label{def:pendant-path}
A \emph{pendant path} in a tree $H$ is a path of the form
\begin{equation}\label{eq:pendant-path}
x-u-v,\quad \text{with} \quad d_H(x)\ge3,\quad d_H(u)=2,\quad \text{and} \quad d_H(v)=1.
\end{equation}
\end{definition}
In Definition \ref{def:pendant-path}, the leaf $v$ is called a \emph{terminal leaf} and a leaf that is not a terminal leaf is called a \emph{direct leaf}.

The next observation gives a restriction on how branch vertices can be connected in the component \(T\).

\begin{lemma}\label{lem:deg3-no-two-branch-neighbors}
A degree-\(3\) branch vertex of the extremal component \(T\) cannot have two branch neighbors.
\end{lemma}

\begin{figure}[ht]
\centering
\begin{tikzpicture}[
    scale=0.85,
    vertex/.style={circle,draw=black,minimum size=8mm,inner sep=0pt}
]
\node[vertex] (v) at (0,0) {$v$};
\node[vertex] (x) at (-2,0) {$x$};
\node[vertex] (y) at (2,0) {$y$};
\node[vertex] (z) at (0,-1.8) {$z$};
\draw[line width=0.45mm] (v)--(x) (v)--(y) (v)--(z);
\draw[dashed,line width=0.4mm] (x)--(-3.2,0.9) (x)--(-3.2,-0.9);
\draw[dashed,line width=0.4mm] (y)--(3.2,0.9) (y)--(3.2,-0.9);
\node at (-2,1.45) {$d_T(x)\ge3$};
\node at (2,1.45) {$d_T(y)\ge3$};
\end{tikzpicture}
\caption{A forbidden configuration of branch vertices in an extremal component.}
\label{fig:forbidden-branch}
\end{figure}
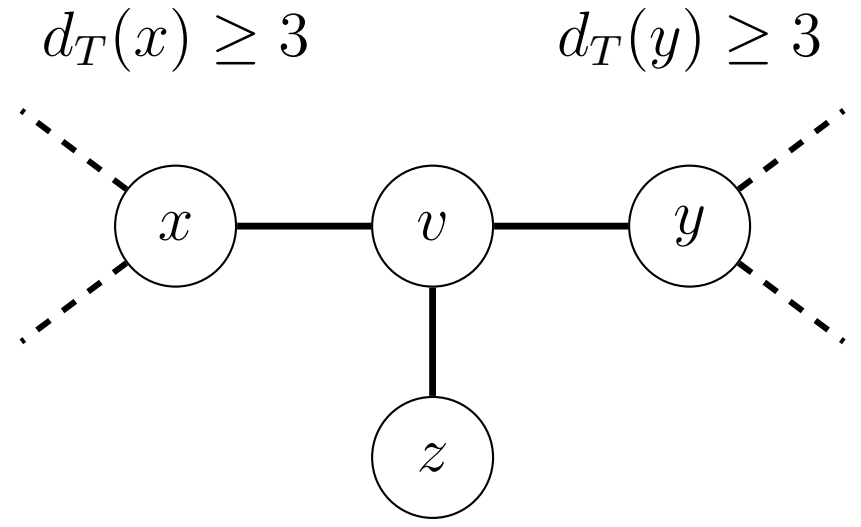

\begin{proof}
Suppose that \(d_T(v)=3\) and that \(x,y\) are two branch neighbors of \(v\) as shown in Figure~\ref{fig:forbidden-branch}. Since \(T\) is a tree, \((x,y)\notin E(T)\). We apply the involution \((v,x,y)\). The component \(T=\Gamma_d\) contains \((v,x)\) and \((v,y)\), so one of these edges is replaced by \((x,y)\). Thus \(v\) becomes a degree-\(2\) vertex. One of \(x\) or \(y\) gains one degree and the other retains its degree, so both remain branch vertices. No leaf is created or lost. The resulting \(d\)-th component therefore has the same number of leaves as \(T\) but one fewer branch vertex, contradicting \textup{(E2)}.
\end{proof}

Lemma \ref{lem:deg3-no-two-branch-neighbors} shows that the component \(T\) can not have the configuration illustrated in Figure~\ref{fig:forbidden-branch}.

We now apply Corollary~\ref{cor:dminusone} to the terminal
leaf of a pendant path in $T$ and use the extremal
conditions to restrict its possible destinations.

\begin{lemma}\label{lem:direct-leaves}
If the extremal component \(T\) satisfies \(p(T)\ge1\), then
\[
\ell(T)-p(T)\ge d-1.
\]
\end{lemma}

\begin{proof}
We choose a pendant path \(x-u-v\) in \(T\) as in \eqref{eq:pendant-path}. By Corollary~\ref{cor:dminusone}, the terminal leaf \(v\) has at least \(d-1\) possible destinations.

A destination cannot be a non-leaf. Indeed, if the involution replaces \((u,v)\) by \((w,v)\), where \(w\) is a non-leaf, then \(u\) becomes a new leaf, \(v\) remains a leaf, and no old leaf is lost. This produces a component with more leaves than \(T\), contradicting \textup{(E1)}.

A destination also cannot be the terminal leaf \(w\) of another pendant path \(y-z-w\). In that case, replacing \((u,v)\) by \((w,v)\) makes \(u\) a leaf and \(w\) a degree-\(2\) vertex. Thus the numbers of leaves and branch vertices remain unchanged, so the resulting component is again extremal. However, \(w\) and \(z\) become adjacent degree-\(2\) vertices, contradicting Corollary~\ref{cor:no-adjacent-deg2}.

It follows that every successful destination is a leaf adjacent to a branch vertex. By Lemma~\ref{lem:deg2-adj-leaf} and Corollary~\ref{cor:no-adjacent-deg2}, each degree-\(2\) vertex is adjacent to exactly one leaf. Therefore, the number of direct leaves is
$\ell(T)-p(T)$. Moreover, since every destination of $v$ must be a direct leaf
of $T$, and there exists at least
$d-1$ distinct destinations, we obtain
\[
\ell(T)-p(T)\ge d-1,
\]
as required.
\end{proof}

The proof of Lemma~\ref{lem:direct-leaves} also yields
the following corollary.

\begin{corollary}\label{cor:direct-destinations}
The extremal component \(T\) satisfies that; if $v$ is
the terminal leaf of a pendant path in $T$,
then every destination of $v$ is a direct leaf of $T$.
\end{corollary}

We can now establish the main result of this section.

\begin{theorem}\label{thm:general-lower-bound}
The extremal component \(T\) satisfies
\(
\ell(T)\ge d+1.
\)
\end{theorem}

\begin{proof}
Suppose first that \(p(T)=0\). By \eqref{eq:count-lpb}, we have
\[
b(T)=2d-\ell(T).
\]
Since each branch vertex contributes at least \(1\) in the identity of Lemma~\ref{lem:leaf-identity}, we have
\[
\ell(T)\ge2+b(T)=2+2d-\ell(T).
\]
Hence \(\ell(T)\ge d+1\).

Suppose now that \(p(T)\ge1\). By Lemma~\ref{lem:direct-leaves},
\[
\ell(T)-p(T)\ge d-1.
\]
Using \eqref{eq:count-lpb}, this is equivalent to
\[
2\ell(T)+b(T)\ge3d-1.
\]
On the other hand, Lemma~\ref{lem:leaf-identity} gives \(b(T)\le\ell(T)-2\). Combining these inequalities, we obtain
\[
3\ell(T)-2\ge3d-1,
\]
and therefore
\[
3\ell(T)\ge3d+1.
\]
Since \(\ell(T)\) is an integer, it follows that \(\ell(T)\ge d+1\).
\end{proof}

Theorem~\ref{thm:general-lower-bound}, together with Lemma~\ref{lem:twinstar-characterization}, shows that the extremal component \(T\) satisfies
\[
d+1\le\ell(T)\le2d-2.
\]
For \(d=5\), this leaves only three possibilities, which are \(
\ell(T)\in\{6,7,8\}.
\)
By Lemma~\ref{lem:twinstar-characterization}, the case \(\ell(T)=8\) gives that \(T\) is isomorphic to \(TS_5\). Thus to establish
Conjecture~\ref{conj:TSd} for $d=5$, it suffices to exclude the cases
$\ell(T)=6$ and $\ell(T)=7$. We consider these cases in the
following sections.

\section{Preliminary reductions for \texorpdfstring{$d=5$}{d=5}}\label{sec:d5-prelim}

We now specialize to cycle-free $5$-partitions of $K_{10}$. Before examining the six- and seven-leaf cases separately, we establish two results that are useful in both. Throughout this section, we denote by $\mathcal P$ a cycle-free $5$-partition of $K_{10}$ and choose a weakly equivalent partition $(\Gamma_1,\ldots,\Gamma_5)$ whose component $T=\Gamma_5$ satisfies the extremal conditions \textup{(E1)} and \textup{(E2)} of Definition~\ref{def:extremal-component}.

We first determine the possible degrees of the branch vertices in the two cases $\ell(T)=6$ and $\ell(T)=7$. 

\begin{lemma}\label{lem:d5-shapes}
If $\ell(T)=6$, the possible branch-degree multisets are of type
\[
\{5,3\},\quad \{4,4\},\quad \text{or} \quad \{4,3,3\}.
\]
If $\ell(T)=7$, the possible branch-degree multisets are of type
\[
\{5,4\},\quad \{5,3,3\},\quad \text{or} \quad \{4,4,3\}.
\]
\end{lemma}

\begin{proof}
Suppose first that $\ell(T)=6$. Lemma~\ref{lem:leaf-identity} gives
\[
\sum_{d_T(v)\ge3}\bigl(d_T(v)-2\bigr)=4.
\]
By Lemma~\ref{lem:degree-bound}, every branch degree is at most $5$, so the summands belong to $\{1,2,3\}$. The possible decompositions are
\[
3+1,\qquad 2+2,\qquad 2+1+1,\qquad 1+1+1+1,
\]
corresponding to $\{5,3\}$, $\{4,4\}$, $\{4,3,3\}$, and $\{3,3,3,3\}$ branch-degree multisets, respectively. The last case is impossible; indeed, this case has the property that $\ell(T)=6$ and $b(T)=4$ which implies that $p(T)=0$ by \eqref{eq:count-lpb}. Therefore, the four branch vertices induce a tree. One of the branch vertices must have at least two branch neighbors, contradicting Lemma~\ref{lem:deg3-no-two-branch-neighbors}.

Suppose now that $\ell(T)=7$. Since $T$ has only three non-leaf vertices, $b(T)\le3$, and Lemma~\ref{lem:leaf-identity} yields
\[
\sum_{d_T(v)\ge3}\bigl(d_T(v)-2\bigr)=5.
\]
The only decompositions of $5$ into at most three summands from $\{1,2,3\}$ are $3+2$, $3+1+1$, and $2+2+1$. These give the stated multisets.
\end{proof}

The following result is used to exclude the configurations $\{4,3,3\}$ and $\{5,3,3\}$ in the six- and seven-leaf cases, respectively.

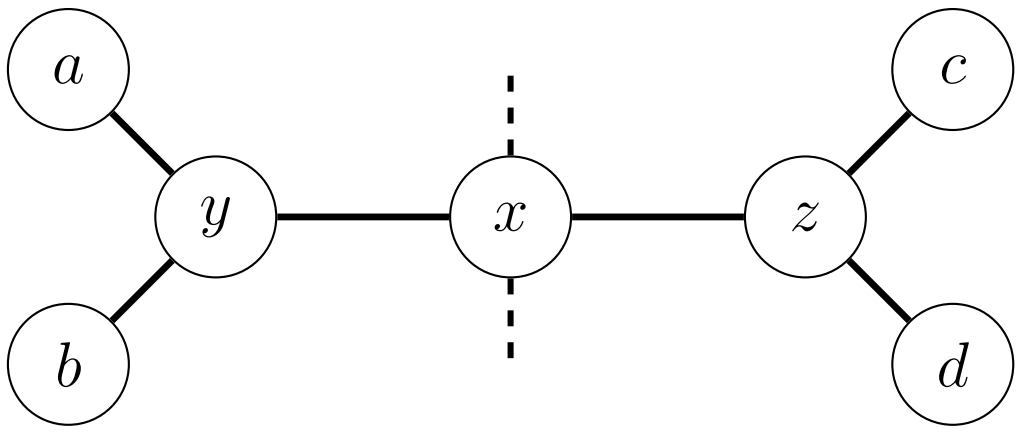
\begin{figure}[ht]
\centering
\begin{tikzpicture}[
    scale=0.65,
    vertex/.style={circle,draw=black,minimum size=8mm,inner sep=0pt}
]
\node[vertex] (y) at (-3,0) {$y$};
\node[vertex] (x) at (0,0) {$x$};
\node[vertex] (z) at (3,0) {$z$};
\node[vertex] (a) at (-4.5,1.5) {$a$};
\node[vertex] (b) at (-4.5,-1.5) {$b$};
\node[vertex] (c) at (4.5,1.5) {$c$};
\node[vertex] (d) at (4.5,-1.5) {$d$};
\foreach \s/\t in {y/x,x/z,y/a,y/b,z/c,z/d}
    \draw[line width=0.45mm] (\s)--(\t);
\draw[dashed,line width=0.4mm] (x)--(0,1.55);
\draw[dashed,line width=0.4mm] (x)--(0,-1.55);
\end{tikzpicture}
\caption{Two nonadjacent degree-$3$ branch vertices with two leaf neighbors each and $d_T(x) > 3$.}
\label{fig:d5-five-components}
\end{figure}

\begin{lemma}\label{lem:d5-five-components}
Suppose that $y$ and $z$ are nonadjacent degree-$3$
branch vertices of $T$, with leaf neighbors $a,b$
at $y$ and $c,d$ at $z$, as illustrated in
Figure~\ref{fig:d5-five-components}.
Then such a configuration is impossible.
\end{lemma}

\begin{proof}
Observe first that transferring either $a$ or $b$
from $y$ to $z$, or either $c$ or $d$ from $z$ to $y$,
would preserve the number of leaves and decrease the
number of branch vertices of $T$. Indeed, the vertex
losing the leaf would have degree $2$, while the other
would have degree $4$. Thus none of these transfers
is possible, by \textup{(E2)}.

Since $y$ and $z$ are nonadjacent in $T$, the five edges
\[
(y,z),\quad (z,a),\quad (z,b),\quad
(y,c),\quad \text{and}\quad (y,d)
\]
belong to the other four components. We show that
these five edges must belong to different components.

First, $(y,z)$ cannot belong to the same component
as any of the other four edges. For instance, if
$(y,z)$ and $(z,a)$ belonged to the same component,
then Lemma~\ref{lem:leaf-switch} would imply that
the involution $(y,z,a)$ transfers the leaf $a$
from $y$ to $z$, contradicting \textup{(E2)}.
The same argument applies to the other three edges.

Now suppose that two of the four edges
\[
(z,a),\quad (z,b),\quad (y,c),\quad (y,d)
\]
belong to the same component $\Gamma_r$, where
$r\in\{1,2,3,4\}$. Recall that $T=\Gamma_5$.
By interchanging $y$ and $z$ and relabeling
their leaf neighbors if necessary, we may assume
that $(z,a)\in E(\Gamma_r)$ and that either
$(z,b)$ or $(y,c)$ also belongs to $\Gamma_r$.

We apply the involution $(y,z,a)$ to the original
partition $(\Gamma_1,\ldots,\Gamma_5)$ and denote
the resulting partition by $(\Delta_1,\ldots,\Delta_5)$.
Since transferring the leaf $a$ from $y$ to $z$
would contradict \textup{(E2)},
Lemma~\ref{lem:leaf-switch} implies that the
involution leaves the fifth component $T$ unchanged
and exchanges the component memberships of
$(y,z)$ and $(z,a)$. In particular,
\[
\Delta_5=T,\quad \text{and} \quad (y,z)\in E(\Delta_r).
\]

If $(z,b)\in E(\Gamma_r)$, then
$(z,b)\in E(\Delta_r)$ as well, since the
involution $(y,z,a)$ does not involve this edge.
Thus $(y,z)$ and $(z,b)$ belong to the same
component $\Delta_r$. By Lemma~\ref{lem:leaf-switch},
the involution $(y,z,b)$ transfers the leaf $b$
from $y$ to $z$ in the fifth component $T$,
contradicting \textup{(E2)}.

Similarly, if $(y,c)\in E(\Gamma_r)$,
then $(y,c)\in E(\Delta_r)$. Since
$(y,z)$ also belongs to $\Delta_r$,
the involution $(z,y,c)$ transfers the leaf
$c$ from $z$ to $y$ in $T$, again
contradicting \textup{(E2)}.

Therefore, the four edges must belong to
different components. Together with the
first part of the proof, this shows that
all five edges belong to different components
among $\Gamma_1,\ldots,\Gamma_4$, which is
impossible.
\end{proof}

In the next section, we turn to the particular case where $\ell(T)=6$.

\section{The six-leaf case}\label{sec:d5-six}

In this section, we still consider $\mathcal P$ to be a cycle-free $5$-partition of $K_{10}$ and choose a weakly equivalent partition $(\Gamma_1,\ldots,\Gamma_5)$ whose component $T=\Gamma_5$ satisfies the extremal conditions \textup{(E1)} and \textup{(E2)}. We focus on the case where $\ell(T)=6$ and show that this is not possible. By Lemma~\ref{lem:d5-shapes}, there are only three possible branch-degree multisets: $\{5,3\}$, $\{4,3,3\}$, and $\{4,4\}$. We start with the following lemma.

\begin{lemma}\label{lem:d5-six-mobility}
Suppose $\ell(T)=6$ and $p(T)=2$. Then the terminal leaf of either pendant path can be transferred to every leaf adjacent to a branch vertex.
\end{lemma}

\begin{proof}
There are exactly $\ell(T)-p(T)=4$ leaves adjacent to branch vertices. By Corollary~\ref{cor:dminusone}, the terminal leaf of either pendant path has at least $d-1=4$ destinations. From Corollary \ref{cor:direct-destinations}, neither a non-leaf nor the terminal leaf of the other pendant path can be a destination. Thus all four leaves adjacent to branch vertices are destinations.
\end{proof}

\subsection{The branch-degree multiset \texorpdfstring{$\{5,3\}$}{{5,3}}}

\begin{lemma}\label{lem:d5-no53}
If the extremal component \(T\) satisfies \(\ell(T)= 6\), then \(T\) cannot have branch-degree multiset $\{5,3\}$.
\end{lemma}

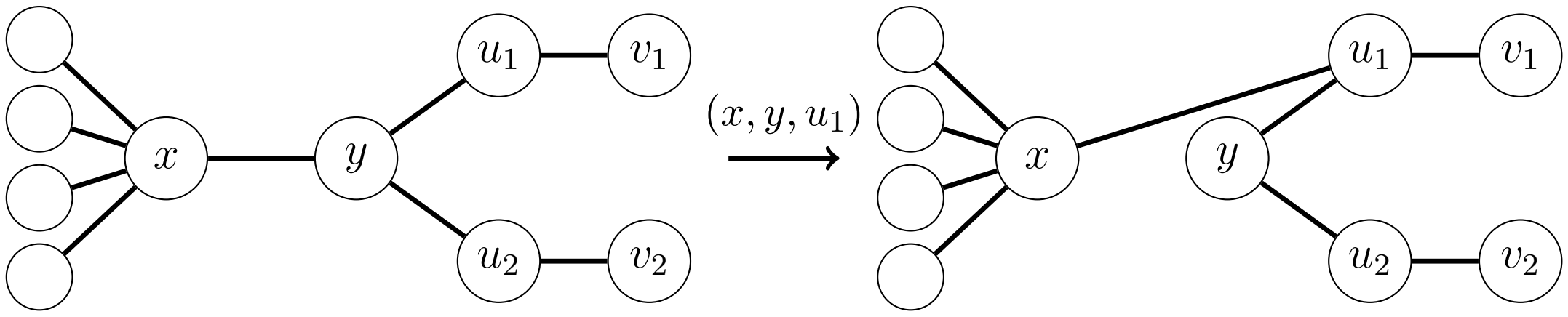
\begin{figure}[ht]
\centering
\begin{tikzpicture}[
    scale=0.72,
    vertex/.style={circle,draw=black,minimum size=7.5mm,inner sep=0pt},
    leaf/.style={circle,draw=black,minimum size=6mm,inner sep=0pt}
]

\begin{scope}
\node[vertex] (x) at (-2.4,0) {$x$};
\node[vertex] (y) at (0,0) {$y$};
\node[vertex] (u1) at (1.8,1.3) {$u_1$};
\node[vertex] (v1) at (3.7,1.3) {$v_1$};
\node[vertex] (u2) at (1.8,-1.3) {$u_2$};
\node[vertex] (v2) at (3.7,-1.3) {$v_2$};

\foreach \p/\q in {x/y,y/u1,u1/v1,y/u2,u2/v2}
    \draw[line width=0.45mm] (\p)--(\q);

\foreach \h [count=\i] in {-1.5,-0.5,0.5,1.5} {
    \node[leaf] (L\i) at (-4.0,\h) {};
    \draw[line width=0.45mm] (x)--(L\i);
}
\end{scope}

\draw[->,line width=0.45mm] (4.7,0)--(6.1,0);
\node at (5.4,0.55) {\small $(x,y,u_1)$};

\begin{scope}[shift={(11,0)}]
\node[vertex] (x) at (-2.4,0) {$x$};
\node[vertex] (y) at (0,0) {$y$};
\node[vertex] (u1) at (1.8,1.3) {$u_1$};
\node[vertex] (v1) at (3.7,1.3) {$v_1$};
\node[vertex] (u2) at (1.8,-1.3) {$u_2$};
\node[vertex] (v2) at (3.7,-1.3) {$v_2$};

\foreach \p/\q in {x/u1,y/u1,u1/v1,y/u2,u2/v2}
    \draw[line width=0.45mm] (\p)--(\q);

\foreach \h [count=\i] in {-1.5,-0.5,0.5,1.5} {
    \node[leaf] (R\i) at (-4.0,\h) {};
    \draw[line width=0.45mm] (x)--(R\i);
}
\end{scope}

\end{tikzpicture}
\caption{The involution $(x,y,u_1)$ creates a degree-$2$
vertex $y$ with no leaf neighbor.}
\label{fig:d5-no53-switch}
\end{figure}

\begin{proof}
Suppose that $x$ and $y$ are the two branch vertices, with $d_T(x)=5$ and $d_T(y)=3$, as shown in Figure~\ref{fig:d5-no53-switch}. They must be adjacent: otherwise the path joining them would contain a degree-$2$ vertex adjacent to no leaf, contrary to Lemma~\ref{lem:deg2-adj-leaf}. Since $\ell(T)=6$ and $b(T)=2$, equation~\eqref{eq:count-lpb} gives $p(T)=2$. Hence $T$ has exactly two pendant paths.

By Lemma~\ref{lem:d5-six-mobility}, we may successively relocate the two pendant paths so that both are based at $y$. After relabeling, the two paths have the form
\[
y-u_1-v_1,\quad \text{and} \quad y-u_2-v_2.
\]

We apply the triangle involution $(x,y,u_1)$. The component $T$ contains $(x,y)$ and $(y,u_1)$, but not $(x,u_1)$. Replacing $(y,u_1)$ by $(x,u_1)$ would raise the degree of $x$ from $5$ to $6$, contradicting Lemma~\ref{lem:degree-bound}. Hence the involution replaces $(x,y)$ by $(x,u_1)$. In the resulting fifth component, $y$ has degree $2$ and is adjacent to the non-leaves $u_1$ and $u_2$. As shown in Figure~\ref{fig:d5-no53-switch}, the component still has six leaves and two branch vertices, so it remains extremal, but violates Lemma~\ref{lem:deg2-adj-leaf}.
\end{proof}

\subsection{The branch-degree multiset \texorpdfstring{$\{4,3,3\}$}{{4,3,3}}}

\begin{lemma}\label{lem:d5-no433-six}
If the extremal component \(T\) satisfies \(\ell(T) = 6\), then \(T\) cannot have branch-degree multiset $\{4,3,3\}$.
\end{lemma}

\begin{figure}[ht]
\centering
\begin{tikzpicture}[
    scale=0.65,
    vertex/.style={circle,draw=black,minimum size=7.5mm,inner sep=0pt},
    leaf/.style={circle,draw=black,minimum size=6mm,inner sep=0pt}
]

\node[vertex] (x) at (0,0) {$x$};
\node[vertex] (y) at (-2.6,0) {$y$};
\node[vertex] (z) at (2.6,0) {$z$};

\node[vertex] (u) at (0,-1.8) {$u$};
\node[vertex] (v) at (0,-3.6) {$v$};

\node[vertex] (a) at (0,2.0) {$a$};

\node[vertex] (b) at (-4.5,1.3) {$b$};
\node[vertex] (c) at (-4.5,-1.3) {$c$};

\node[vertex] (d) at (4.5,1.3) {$d$};
\node[vertex] (e) at (4.5,-1.3) {$e$};

\foreach \p/\q in {y/x,x/z,x/a,x/u,u/v,y/b,y/c,z/d,z/e}
    \draw[line width=0.45mm] (\p)--(\q);

\end{tikzpicture}
\caption{The extremal component \(T\) with branch-degree multiset
\(\{4,3,3\}\).}
\label{fig:d5-433-pendant-at-x}
\end{figure}
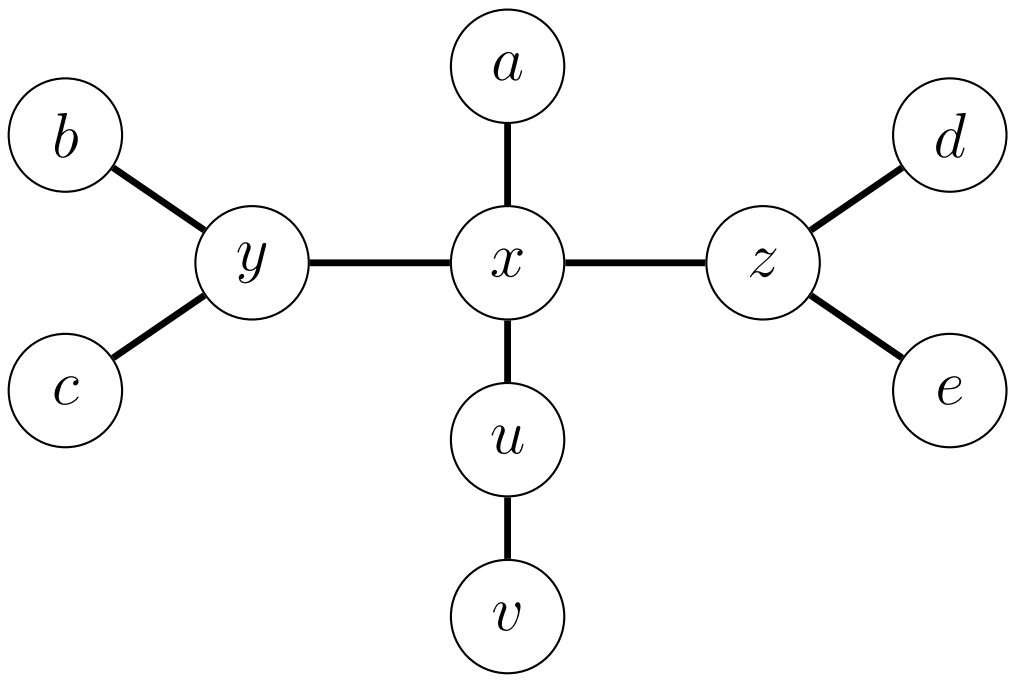

\begin{proof}
Here $b(T)=3$, so \eqref{eq:count-lpb} gives $p(T)=1$. By Lemma~\ref{lem:deg3-no-two-branch-neighbors}, the degree-$4$ vertex must be the middle vertex. We may therefore write
\[
y-x-z,\qquad d_T(x)=4,\quad d_T(y)=d_T(z)=3.
\]

Since $p(T)=1$, there is a unique pendant path. 

We show that, after passing to a weakly equivalent partition,
we may assume that the unique pendant path is based at $x$, as illustrated
in Figure~\ref{fig:d5-433-pendant-at-x}.

If the pendant path is already based at $x$, there is
nothing to prove. Otherwise, after interchanging $y$ and $z$ if necessary, we can write it as $y-b-v$, where $v$ is the terminal leaf. The remaining five leaves are adjacent directly to branch vertices. By Corollary~\ref{cor:dminusone} and \ref{cor:direct-destinations}, at least four of these five leaves receive $v$. Since $x$ has two direct leaf neighbors, at least one is a destination. Transferring $v$ to this leaf moves the pendant path
from $y$ to $x$, without changing the numbers of leaves
or branch vertices. Thus the resulting component remains
extremal and has the configuration shown in
Figure~\ref{fig:d5-433-pendant-at-x}.

We may thus assume that the unique pendant path is based at $x$. Consequently, each of $y$ and $z$ has exactly two leaf neighbors. Lemma~\ref{lem:d5-five-components} excludes this configuration.
\end{proof}

\subsection{The branch-degree multiset
\texorpdfstring{$\{4,4\}$}{{4,4}}}

The remaining case requires us to track the components containing certain edges outside $T$. To do this, we first choose a convenient form for the component \(T\) of \(\mathcal{P}\).

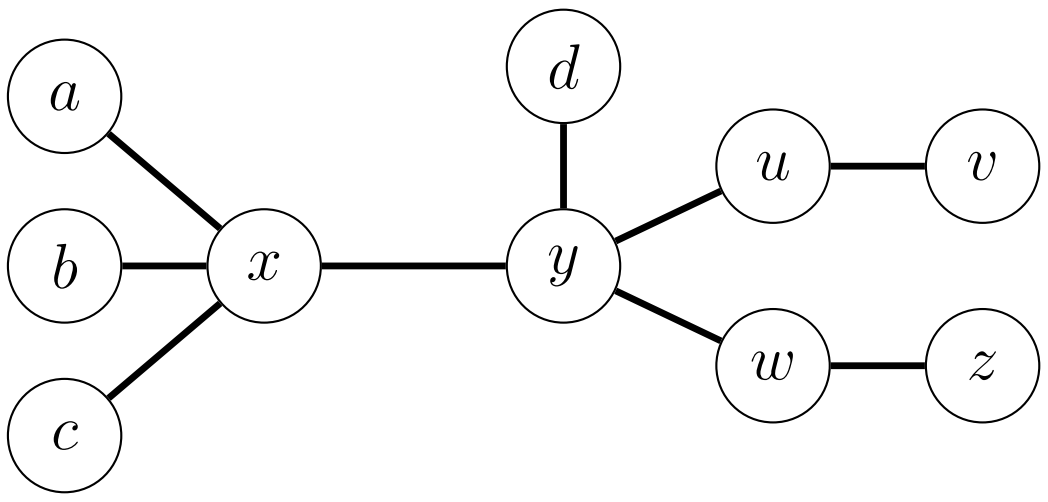
\begin{figure}[ht]
\centering
\begin{tikzpicture}[
    scale=0.66,
    vertex/.style={circle,draw=black,minimum size=7.5mm,inner sep=0pt}
]
\node[vertex] (x) at (0,0) {$x$};
\node[vertex] (y) at (3,0) {$y$};
\node[vertex] (a) at (-2,1.7) {$a$};
\node[vertex] (b) at (-2,0) {$b$};
\node[vertex] (c) at (-2,-1.7) {$c$};
\node[vertex] (d) at (3,2) {$d$};
\node[vertex] (u) at (5.1,1) {$u$};
\node[vertex] (v) at (7.2,1) {$v$};
\node[vertex] (w) at (5.1,-1) {$w$};
\node[vertex] (z) at (7.2,-1) {$z$};
\foreach \p/\q in {x/y,x/a,x/b,x/c,y/d,y/u,u/v,y/w,w/z}
    \draw[line width=0.45mm] (\p)--(\q);
\end{tikzpicture}
\caption{A convenient form of the extremal six-leaf component \(T\) with branch-degree multiset $\{4,4\}$.}
\label{fig:d5-normal44}
\end{figure}

\begin{lemma}\label{lem:d5-normal44}
If the extremal component \(T\) satisfies $\ell(T)=6$ and the branch-degree multiset of \(T\) is $\{4,4\}$, then after passing to a weakly equivalent partition, we may assume that
\begin{equation}\label{eq:d5-normal44}
\begin{aligned}
E(T)=\{(x,y),(x,a),(x,b),(x,c),(y,d),
       (y,u),(u,v),(y,w),(w,z)\}.
\end{aligned}
\end{equation}
Thus $x,y$ have degree $4$, the vertices $u,w$ have degree $2$, and $a,b,c,d,v,z$ are leaves, as shown in Figure~\ref{fig:d5-normal44}.
\end{lemma}

\begin{proof}
Since $\ell(T)=6$ and $b(T)=2$, equation~\eqref{eq:count-lpb} gives $p(T)=2$. The two branch vertices are adjacent, because no degree-$2$ vertex lies between two branch vertices. By Lemma~\ref{lem:d5-six-mobility}, we can successively relocate the pendant paths so that both are based at the same branch vertex, say $y$. At each step, the number of leaves and branch vertices is preserved, so the resulting fifth component remains extremal. Relabeling the vertices gives \eqref{eq:d5-normal44}.
\end{proof}

We now consider the component \(T = \Gamma_5\) as in Figure~\ref{fig:d5-normal44}. The next lemma shows that the edges \((u,w)\), \((w,v)\), \((v,z)\), and \((z,u)\) belong to four distinct components and, more importantly, that their component memberships can be permuted while \(T\) remains unchanged.

\begin{lemma}\label{lem:d5-square}
If the extremal component $T$ is as in
Figure~\ref{fig:d5-normal44}, then the four edges
\[
(u,w),\quad (w,v),\quad (v,z),\quad \text{and}\quad (z,u)
\]
belong to four different components among
$\Gamma_1,\ldots,\Gamma_4$. Moreover, for any ordering $i,j,k,l$ of
$1,2,3,4$, there exists a finite sequence of
triangle involutions producing a partition
$(\Delta_1,\ldots,\Delta_5)$ such that
$\Delta_5=T$ and
\[
\begin{aligned}
(u,w)\in E(\Delta_i), \quad
(w,v)\in E(\Delta_j), \quad
(v,z)\in E(\Delta_k), \quad \text{and} \quad
(z,u)&\in E(\Delta_l).
\end{aligned}
\]
All these involutions involve only the vertices
$u,v,w,z$.
\end{lemma}

\begin{proof}
Since none of the four edges belongs to \(T=\Gamma_5\), we write
\[
\begin{aligned}
(u,w)\in E(\Gamma_i), \quad (w,v)\in E(\Gamma_j),\quad
(u,z)\in E(\Gamma_k), \quad \text{and} \quad  (v,z)\in E(\Gamma_l),
\end{aligned}
\]
where \(i,j,k,l\in\{1,2,3,4\}\). We first show that these indices
are pairwise distinct.

Consider the involution \((u,w,v)\). The only edge of this
triangle belonging to \(T\) is \((u,v)\). Replacing it by \((w,v)\)
would create an additional leaf, contradicting \textup{(E1)},
while replacing it by \((u,w)\) would isolate \(v\).
Hence the involution leaves \(T\) unchanged and exchanges
\((u,w)\) and \((w,v)\), which forces \(i\ne j\).

Similarly, the involution \((u,z,v)\) cannot replace \((u,v)\)
by \((u,z)\), since this would isolate \(v\). Replacing
\((u,v)\) by \((z,v)\) would preserve the number of leaves and
branch vertices but create adjacent degree-\(2\) vertices
\(w\) and \(z\), contradicting Corollary~\ref{cor:no-adjacent-deg2}.
Thus \(T\) remains unchanged, and \((u,z)\) and \((v,z)\)
are exchanged, forcing \(k\ne l\).

The same reasoning applied to the involutions \((w,u,z)\)
and \((w,v,z)\) shows that both leave \(T\) unchanged and
exchange \((u,w)\) with \((u,z)\), and \((w,v)\) with
\((v,z)\), respectively. Consequently,
\[
i\ne k,\quad \text{and} \quad j\ne l.
\]

We now prove that $j\ne k$ and $i\ne l$.
First, apply the involution $(u,w,v)$ to the original
partition and denote the result by
$(\Delta_1,\ldots,\Delta_5)$. This involution
leaves $T$ unchanged and exchanges the component
memberships of $(u,w)$ and $(w,v)$. Hence
\[
\Delta_5=T,\qquad
(u,w)\in E(\Delta_j),\quad \text{and} \quad
(u,z)\in E(\Delta_k).
\]
We then apply the involution $(w,u,z)$ to this
new partition. Since it must also leave $T$
unchanged, it exchanges the component memberships
of $(u,w)$ and $(u,z)$. These two edges must
therefore belong to different components, giving
$j\ne k$.

To show that $i\ne l$, we return to the original
partition and apply the involution $(u,z,v)$.
Denote the resulting partition by
$(\Theta_1,\ldots,\Theta_5)$. This involution
leaves $T$ unchanged and exchanges the component
memberships of $(u,z)$ and $(v,z)$, so
\[
\Theta_5=T,\qquad
(u,z)\in E(\Theta_l),\quad \text{and} \quad
(u,w)\in E(\Theta_i).
\]
Applying the involution $(w,u,z)$ to this partition,
we again obtain an exchange of the component
memberships of $(u,w)$ and $(u,z)$ while
leaving $T$ unchanged. Consequently, these
edges must belong to different components,
and hence $i\ne l$.

Therefore \(i,j,k,l\) are pairwise distinct, and hence
\(
\{i,j,k,l\}=\{1,2,3,4\}.
\)

Finally, the involutions \((u,w,v)\), \((u,z,v)\),
\((w,u,z)\), and \((w,v,z)\) exchange the component
memberships of the following pairs of edges respectively:
\[
\begin{gathered}
((u,w)\ (w,v)),\quad ((u,z)\ (v,z)), \quad
((u,w)\ (u,z)),\quad  \text{and} \quad ((w,v)\ (v,z)).
\end{gathered}
\]
These exchanges generate the symmetric group $S_4$,
so every permutation of the component memberships
of the four edges can be realized. Since each
involution leaves $T$ unchanged, the same holds
for any finite sequence of these involutions.
\end{proof}

A consequence of Lemma \ref{lem:d5-square} is the following result.

\begin{corollary} \label{lem:d5-xedges}
If the extremal component $T$ is as in
Figure~\ref{fig:d5-normal44}, then the edges
\[
(x,u),\quad(x,v),\quad(x,w),\quad \text{and} \quad (x,z)
\]
belong to four different components among $\Gamma_1,\ldots,\Gamma_4$.
\end{corollary}

\begin{proof}
Suppose first that $(x,u)$ and $(x,v)$ belong to the same component.
By Lemma~\ref{lem:leaf-switch}, the involution $(u,x,v)$ replaces
$(u,v)$ by $(x,v)$ in $T$. This makes $u$ a new leaf while $v$
remains a leaf, contradicting \textup{(E1)}. Hence $(x,u)$ and
$(x,v)$ belong to different components. Similarly, $(x,w)$ and
$(x,z)$ belong to different components.

Now suppose that $(x,r)$ and $(x,s)$ belong to the same component
$\Gamma_i$, where $r\in\{u,v\}$ and $s\in\{w,z\}$.
The edge $(r,s)$ is one of the four edges considered in
Lemma~\ref{lem:d5-square}. By that lemma, we can apply a sequence
of involutions to obtain a partition
$(\Delta_1,\ldots,\Delta_5)$ such that $\Delta_5=T$ and
$(r,s)\in E(\Delta_i)$. Since none of these involutions involves
$x$, the edges $(x,r)$ and $(x,s)$ remain in the $i$-th component.
Consequently, the three edges $(x,r)$, $(r,s)$, and $(s,x)$
form a cycle in $\Delta_i$, contradicting the fact that
$\Delta_i$ is a tree. Therefore all four edges belong to
different components.
\end{proof}

We are now ready to exclude the six-leaf branch-degree multiset $\{4,4\}$. 

\begin{lemma}\label{lem:d5-no44-six}
If the extremal component \(T\) satisfies \(\ell(T) =6\), then \(T\) cannot have branch-degree multiset $\{4,4\}$.
\end{lemma}

\begin{proof}
We recall that $T$ is as in Figure~\ref{fig:d5-normal44}. By Corollary~\ref{lem:d5-xedges}, the edges $(x,u)$, $(x,v)$, $(x,w)$, and $(x,z)$ belong to four different components among $\Gamma_1,\ldots,\Gamma_4$. The five edges incident with $x$ but not contained in $T$ are exactly
\[
(x,d),\quad(x,u),\quad(x,v),\quad(x,w),\quad \text{and} \quad (x,z).
\]
Consequently, there is a unique $r\in\{u,v,w,z\}$ such that $(x,d)$ and $(x,r)$ lie in the same component.

We show that, after passing to a weakly equivalent
partition, we may assume that $y-d-r$ is a pendant
path in $T$, while $(x,d)$ and $(x,r)$ remain in
the same component.

Suppose first that $r\in\{v,z\}$. If $r=v$,
Lemma~\ref{lem:d5-six-mobility} allows us to transfer
the terminal leaf $v$ of the pendant path $y-u-v$
to the direct leaf $d$. This produces the pendant
path $y-d-v$, as desired. The case $r=z$ is identical.

Suppose now that $r\in\{u,w\}$. If $r=u$, we first
apply the involution $(y,u,v)$. The component $T$
contains $(y,u)$ and $(u,v)$, but not $(y,v)$.
Replacing $(u,v)$ by $(y,v)$ would increase the
number of leaves to seven, contradicting
\textup{(E1)}. Hence the involution replaces
$(y,u)$ by $(y,v)$, transforming the pendant
path $y-u-v$ into $y-v-u$. The vertex $u$ is
now the terminal leaf, and
Lemma~\ref{lem:d5-six-mobility} allows us to
transfer it to $d$, producing the pendant path
$y-d-u$. The case $r=w$ follows in the same way.

In each case, the involutions preserve the number
of leaves and branch vertices, so the resulting
fifth component remains extremal. Moreover, none
of these involutions involves $x$. Consequently,
the edges $(x,d)$ and $(x,r)$ retain their original
component memberships and therefore still belong
to the same component.

We have thus obtained a weakly equivalent partition
whose fifth component contains the pendant path
$y-d-r$, with $(x,d)$ and $(x,r)$ belonging to
the same component. By Lemma~\ref{lem:leaf-switch},
the involution $(d,x,r)$ transfers the terminal
leaf $r$ from $d$ to $x$. As a result, $d$ becomes
a new leaf while $r$ remains a leaf, increasing
the number of leaves of the fifth component to
seven. This contradicts \textup{(E1)} and completes
the proof.
\end{proof}

Therefore, we get the following result.

\begin{proposition}\label{prop:d5-seven}
The extremal component \(T\) satisfies $\ell(T)\ge7$.
\end{proposition}

In the next section, we show that the case $\ell(T)=7$ is not possible.

\section{The seven-leaf case}\label{sec:d5-seven}

Here, we still consider $\mathcal P$ to be a cycle-free $5$-partition of $K_{10}$ and choose a weakly equivalent partition $(\Gamma_1,\ldots,\Gamma_5)$ whose component $T=\Gamma_5$ satisfies the extremal conditions \textup{(E1)} and \textup{(E2)}. By Proposition~\ref{prop:d5-seven}, the extremal component $T$ has at least seven leaves. If it has eight leaves, then $T$ is isomorphic to $TS_5$ by Lemma~\ref{lem:twinstar-characterization}. We therefore suppose throughout this section that $\ell(T)=7$ and show that this is impossible. By Lemma~\ref{lem:d5-shapes}, the only possible branch-degree multisets are
\[
\{5,3,3\},\quad \{4,4,3\},\quad \text{and} \quad \{5,4\}.
\]
We consider them in this order.

\subsection{The branch-degree multiset \texorpdfstring{$\{5,3,3\}$}{{5,3,3}}}

\begin{lemma}\label{lem:d5-no533-seven}
If the extremal component \(T\) satisfies \(\ell(T) = 7\), then \(T\) cannot have branch-degree multiset $\{5,3,3\}$.
\end{lemma}

\begin{figure}[ht]
\centering
\begin{tikzpicture}[
    scale=0.65,
    vertex/.style={circle,draw=black,minimum size=7.5mm,inner sep=0pt},
    leaf/.style={circle,draw=black,minimum size=6mm,inner sep=0pt}
]

\node[vertex] (x) at (0,0) {$x$};
\node[vertex] (y) at (-2.6,0) {$y$};
\node[vertex] (z) at (2.6,0) {$z$};

\node[vertex] (a) at (0,2.0) {$a$};
\node[vertex] (b) at (-1.3,2.0) {$b$};
\node[vertex] (c) at (1.3,2.0) {$c$};

\node[vertex] (d) at (-4.5,1.3) {$d$};
\node[vertex] (e) at (-4.5,-1.3) {$e$};

\node[vertex] (f) at (4.5,1.3) {$f$};
\node[vertex] (g) at (4.5,-1.3) {$g$};

\foreach \p/\q in {y/x,x/z,x/a,x/b,x/c,y/d,y/e,z/f,z/g}
    \draw[line width=0.45mm] (\p)--(\q);

\end{tikzpicture}
\caption{The extremal component \(T\) in the case \(\ell(T)=7\)
with branch-degree multiset \(\{5,3,3\}\).}
\label{fig:d5-533-shape}
\end{figure}
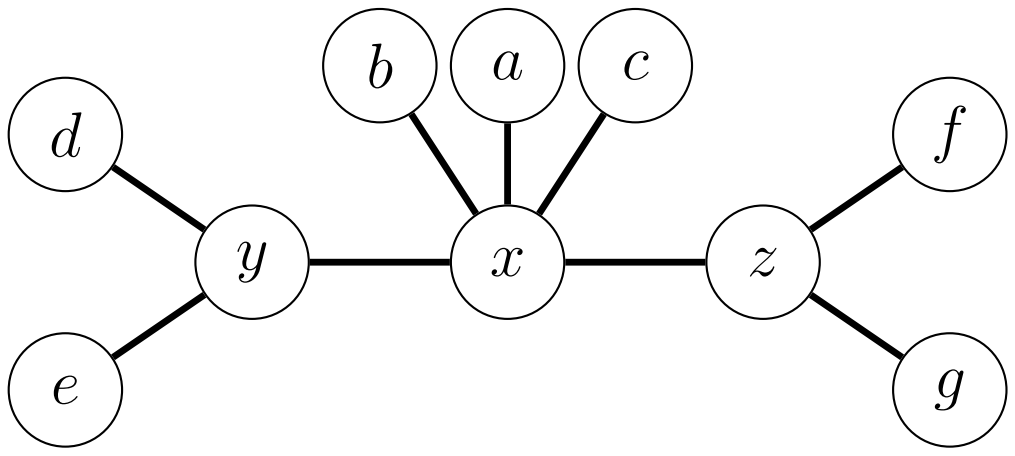

\begin{proof}
Since $\ell(T)=7$ and $b(T)=3$, equation~\eqref{eq:count-lpb}
gives $p(T)=0$. Thus the three branch vertices induce a path.
By Lemma~\ref{lem:deg3-no-two-branch-neighbors}, the middle
vertex of this path cannot have degree $3$. Hence, after relabeling
if necessary, we may write
\[
y-x-z,\qquad d_T(x)=5,\quad d_T(y)=d_T(z)=3.
\]
In particular, each of $y$ and $z$ has exactly two leaf neighbors,
while $x$ has exactly three. This configuration is shown in
Figure~\ref{fig:d5-533-shape}. Lemma~\ref{lem:d5-five-components} shows that this configuration is not possible. 
\end{proof}

\subsection{The branch-degree multiset \texorpdfstring{$\{4,4,3\}$}{{4,4,3}}}

\begin{lemma}\label{lem:d5-no443-seven}
If the extremal component \(T\) satisfies \(\ell(T) = 7\), then \(T\) cannot have branch-degree multiset $\{4,4,3\}$.
\end{lemma}

\begin{figure}[ht]
\centering
\begin{tikzpicture}[
    scale=0.68,
    vertex/.style={circle,draw=black,minimum size=7.5mm,inner sep=0pt}
]
\node[vertex] (x) at (-3,0) {$x$};
\node[vertex] (y) at (0,0) {$y$};
\node[vertex] (z) at (3,0) {$z$};
\node[vertex] (a) at (-4.35,1.7) {$a$};
\node[vertex] (b) at (-5,0) {$b$};
\node[vertex] (c) at (-4.35,-1.7) {$c$};
\node[vertex] (d) at (-0.7,1.85) {$d$};
\node[vertex] (e) at (-0.7,-1.85) {$e$};
\node[vertex] (f) at (4.45,1.35) {$f$};
\node[vertex] (g) at (4.45,-1.35) {$g$};
\foreach \p/\q in {x/y,y/z,x/a,x/b,x/c,y/d,y/e,z/f,z/g}
    \draw[line width=0.45mm] (\p)--(\q);
\end{tikzpicture}
\caption{The seven-leaf configuration with branch degrees $\{4,4,3\}$.}
\label{fig:d5-seven-443}
\end{figure}
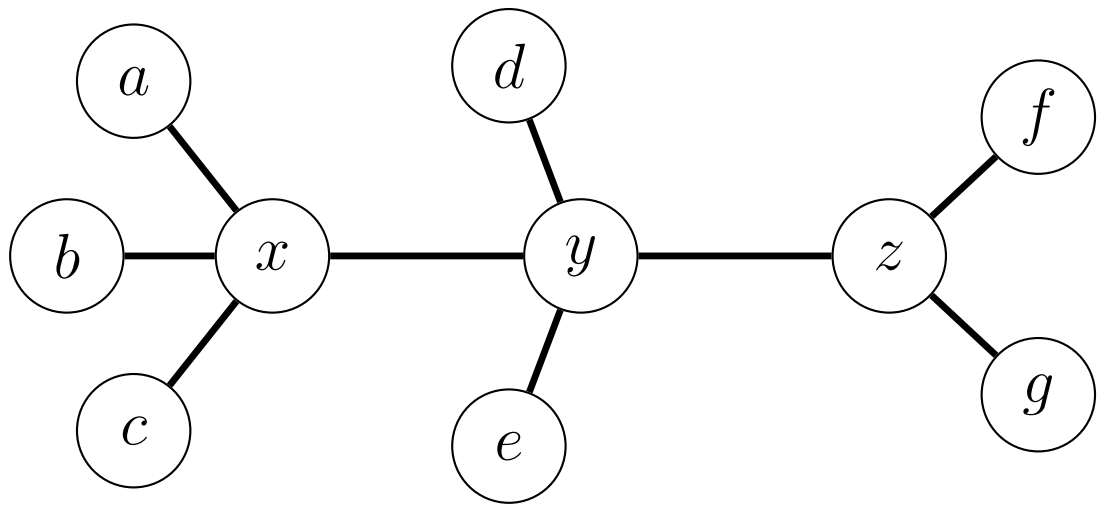

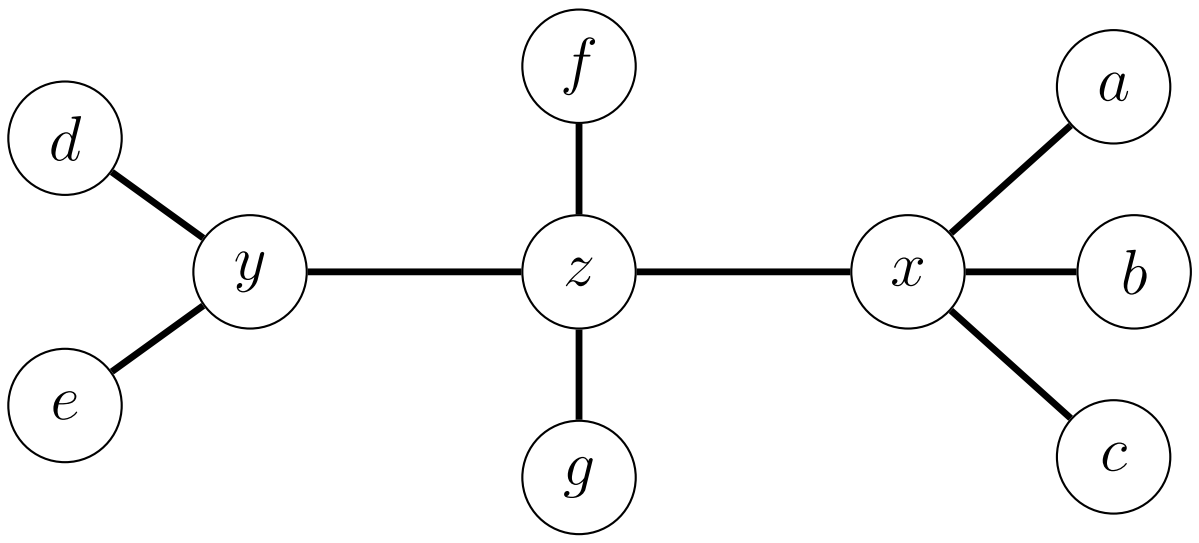
\begin{figure}[ht]
\centering
\begin{tikzpicture}[
    scale=0.68,
    vertex/.style={circle,draw=black,minimum size=7.5mm,inner sep=0pt},
    leaf/.style={circle,draw=black,minimum size=6mm,inner sep=0pt}
]

\node[vertex] (y) at (-3.2,0) {$y$};
\node[vertex] (z) at (0,0) {$z$};
\node[vertex] (x) at (3.2,0) {$x$};

\node[vertex] (d) at (-5.0,1.3) {$d$};
\node[vertex] (e) at (-5.0,-1.3) {$e$};

\node[vertex] (f) at (0,2.0) {$f$};
\node[vertex] (g) at (0,-2.0) {$g$};

\node[vertex] (a) at (5.2,1.8) {$a$};
\node[vertex] (b) at (5.4,0) {$b$};
\node[vertex] (c) at (5.2,-1.8) {$c$};

\foreach \p/\q in {y/z,z/x,y/d,y/e,z/f,z/g,x/a,x/b,x/c}
    \draw[line width=0.45mm] (\p)--(\q);

\end{tikzpicture}
\caption{The component $T'$ with branch path $y-z-x$.}
\label{fig:d5-Tprime-344}
\end{figure}

\begin{proof}
Here, we have that $p(T)=0$, so the three branch vertices form a path. Its middle vertex cannot have degree $3$, by Lemma~\ref{lem:deg3-no-two-branch-neighbors}. Therefore, after relabeling, we write \(T\) as
\[
x-y-z,\qquad d_T(x)=d_T(y)=4,\quad d_T(z)=3,
\]
with
\[
\begin{aligned}
N_T(x)&=\{y,a,b,c\},&N_T(y)&=\{x,z,d,e\},&N_T(z)&=\{y,f,g\}.
\end{aligned}
\]
The configuration is shown in Figure~\ref{fig:d5-seven-443}. The five edges
\[
(x,z),\quad(x,f),\quad(x,g),\quad(x,d),\quad \text{and} \quad (x,e)
\]
do not belong to $T$. We write
\[
\begin{aligned}
(x,z)\in E(\Gamma_i), \quad (x,f)\in E(\Gamma_j), \quad (x,g)\in E(\Gamma_k),\quad
(x,d)\in E(\Gamma_l),\quad  \text{and} \quad (x,e) \in E(\Gamma_m),
\end{aligned}
\]
where $i,j,k,l,m\in\{1,2,3,4\}$. We prove that these indices are pairwise distinct.

First, $i\ne j$ and $i\ne k$. Indeed, if $i=j$,
then $(x,z)$ and $(x,f)$ belong to the same component.
By Lemma~\ref{lem:leaf-switch}, the involution
$(z,x,f)$ transfers the leaf $f$ from $z$ to $x$.
This preserves the number of leaves and decreases
the number of branch vertices, contradicting
\textup{(E2)}. The same argument with $g$ in place
of $f$ gives $i\ne k$.

We now apply the involution $(z,x,f)$ to the original
partition. Since transferring $f$ to $x$ would
contradict \textup{(E2)}, this involution leaves
$T$ unchanged and exchanges the component memberships
of $(x,z)$ and $(x,f)$. Thus $(x,z)$ now belongs
to the $j$-th component, while $(x,g)$ still belongs
to the $k$-th component. If $j=k$, the involution
$(z,x,g)$ would transfer the leaf $g$ from $z$ to $x$,
again contradicting \textup{(E2)}. Hence $j\ne k$.

Moreover, the involutions $(z,x,f)$ and $(z,x,g)$
exchange the component memberships of $(x,z)$
with those of $(x,f)$ and $(x,g)$, respectively,
while leaving $T$ unchanged. These exchanges
correspond to the transpositions $(1,2)$ and
$(1,3)$, which generate the symmetric group $S_3$.
Consequently, starting from the original partition,
we can arrange for $(x,z)$ to belong to any one
of the components $\Gamma_i,\Gamma_j,\Gamma_k$,
without changing $T$ or the component memberships
of $(x,d)$ and $(x,e)$.

Let $t\in\{i,j,k\}$. Choose such a partition
$(\Delta_1,\ldots,\Delta_5)$ with
\[
\Delta_5=T,\qquad
(x,z)\in E(\Delta_t),\qquad
(x,d)\in E(\Delta_l),\qquad
(x,e)\in E(\Delta_m).
\]
We apply the involution $(x,y,z)$. Since replacing
$(y,z)$ by $(x,z)$ would produce a component
with branch-degree multiset $\{5,3,3\}$, which is
impossible by Lemma~\ref{lem:d5-no533-seven},
the involution must replace $(x,y)$ by $(x,z)$.
Denote the resulting partition by
$(\Theta_1,\ldots,\Theta_5)$ and put $T'=\Theta_5$.
Then $T'$ has branch path $y-z-x$, with degrees
$3,4,4$, respectively, as in Figure \ref{fig:d5-Tprime-344} and
\[
(x,y)\in E(\Theta_t),\qquad
(x,d)\in E(\Theta_l),\qquad
(x,e)\in E(\Theta_m).
\]

Transferring either leaf $d$ or $e$ from $y$ to $x$
would preserve the seven leaves and decrease the
number of branch vertices from three to two,
contradicting \textup{(E2)}. Therefore,
Lemma~\ref{lem:leaf-switch} gives
\[
t\ne l,\quad \text{and} \quad t\ne m.
\]
Since $t$ was arbitrary in $\{i,j,k\}$, we obtain
\[
\{i,j,k\}\cap\{l,m\}=\varnothing.
\]

It remains to show that $l\ne m$. Suppose otherwise.
We apply the involution $(y,x,d)$ to the partition
$(\Theta_1,\ldots,\Theta_5)$. Since transferring
$d$ from $y$ to $x$ is forbidden, the involution
leaves $T'$ unchanged and exchanges the component
memberships of $(x,y)$ and $(x,d)$.
Thus $(x,y)$ now belongs to the $l$-th component,
which also contains $(x,e)$ because $l=m$.
By Lemma~\ref{lem:leaf-switch}, the involution
$(y,x,e)$ would then transfer the leaf $e$ from
$y$ to $x$, contradicting \textup{(E2)}.
Hence $l\ne m$.

Hence $i,j,k,l,m$ are five pairwise
distinct elements of $\{1,2,3,4\}$, which is
impossible. \end{proof}

\subsection{The branch-degree multiset \texorpdfstring{$\{5,4\}$}{{5,4}}}

It remains to exclude the case where the branch-degree multiset is $\{5,4\}$. In this case there are two branch vertices and one degree-$2$ vertex. We first choose labels that make it easy to follow the possible transfers of the terminal leaf.

\begin{figure}[ht]
\centering
\begin{tikzpicture}[
    scale=0.62,
    vertex/.style={circle,draw=black,minimum size=7.5mm,inner sep=0pt}
]
\node[vertex] (one) at (-1.4,0) {$1$};
\node[vertex] (two) at (1.6,0) {$2$};
\node[vertex] (four) at (-3.55,2.25) {$4$};
\node[vertex] (six) at (-3.95,0.8) {$6$};
\node[vertex] (eight) at (-3.95,-0.8) {$8$};
\node[vertex] (ten) at (-3.55,-2.25) {$10$};
\node[vertex] (three) at (3.55,1.8) {$3$};
\node[vertex] (five) at (4.0,0) {$5$};
\node[vertex] (seven) at (3.55,-1.8) {$7$};
\node[vertex] (nine) at (-5.5,3.15) {$9$};
\foreach \p/\q in {one/two,one/four,one/six,one/eight,one/ten,two/three,two/five,two/seven}
    \draw[line width=0.44mm] (\p)--(\q);
\draw[line width=0.65mm] (four)--(nine);
\node at (-1.4,-2.85) {};
\end{tikzpicture}
\caption{The tree $H_4$. }
\label{fig:d5-H4}
\end{figure}
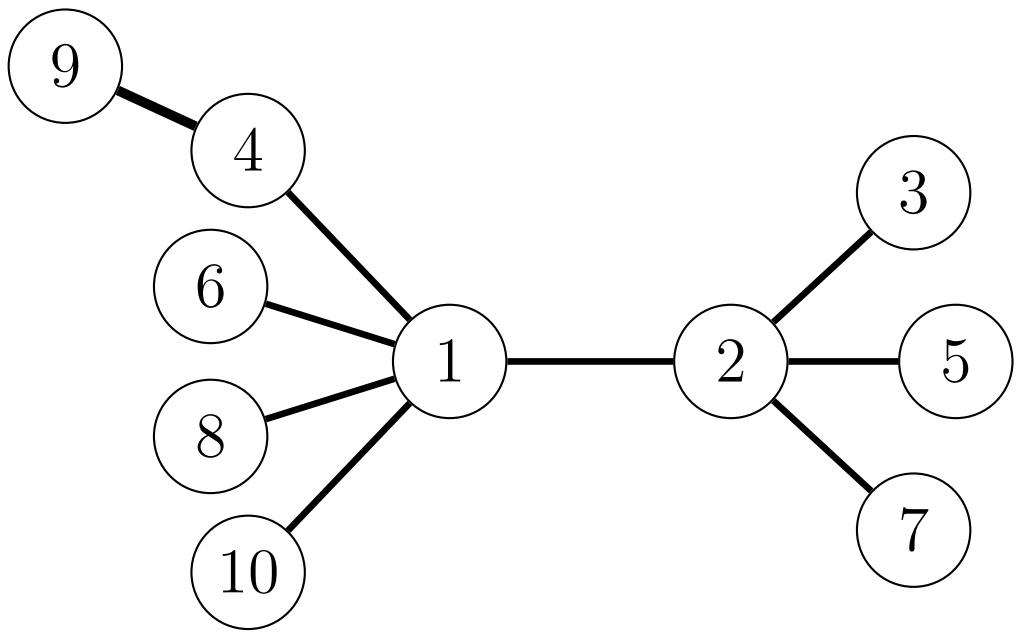

Put
\[
\begin{aligned}
C=\{(1,2),(1,4),(1,6),(1,8),(1,10),
   (2,3),(2,5),(2,7)\},
\end{aligned}
\]
and set
\[
S=\{3,4,5,6,7,8,10\},\quad \text{and} \quad  L=\{4,6,8,10\}.
\]
For $a\in S$, let $H_a$ be the tree defined by
\begin{equation}\label{eq:d5-Ha}
E(H_a)=C\cup\{(a,9)\}.
\end{equation}
Thus $H_a$ consists of the eight fixed edges in $C$ and one edge joining the leaf $9$ to $a$. The vertices in $L$ are adjacent to $1$, whereas $3,5,7$ are adjacent to $2$. Figure~\ref{fig:d5-H4} shows $H_4$ as an example. If the edge $(a,9)$ is replaced by $(2,9)$, the resulting tree is isomorphic to $TS_5$.

\begin{lemma}\label{lem:d5-normal54}
Up to weakly equivalent partition, we may assume that $T=H_a$ for some $a\in L$.
\end{lemma}

\begin{proof}
Let $x,y$ be the branch vertices, with $d_T(x)=5$ and $d_T(y)=4$. By \eqref{eq:count-lpb}, $p(T)=1$, hence the unique degree-$2$ vertex belongs to a pendant path based at $x$ or $y$. If the path is based at $x$, we obtain the required form by relabeling.

Suppose instead that the path is based at $y$. Then $x$ has four direct leaf neighbors and $y$ has two. By Corollary~\ref{cor:dminusone} and \ref{cor:direct-destinations}, the terminal leaf has at least four destinations among these six direct leaves. Therefore at least one direct leaf at $x$ receives it. The transfer moves the pendant path from $y$ to $x$ and preserves the numbers of leaves and branch vertices, so the resulting component is still extremal. Relabeling gives $T=H_a$ for some $a\in L$.
\end{proof}

Our next goal is to show that, while the eight edges of $C$ stay in the fifth component, the leaf $9$ can be attached at several different vertices of $S$. The following two lemmas make this precise.

\begin{lemma}\label{lem:d5-reachable}
Suppose that $\Gamma_5=H_a$ for some $a\in S$, and assume that no cycle-free $5$-partition weakly equivalent to $(\Gamma_1,\ldots,\Gamma_5)$ has a component isomorphic to $TS_5$. Let $R\subseteq S$ be the set of vertices $s$ for which a sequence of triangle involutions, changing only the attachment of the leaf $9$ in the fifth component and keeping all edges of $C$ fixed, produces a partition with fifth component $H_s$. Then $|R|\ge5$.
\end{lemma}

\begin{proof}
Clearly $a\in R$. Suppose that $|R|\le4$, and fix $s\in R$. We consider the weakly equivalent partition whose fifth component is $H_s$. Choose three distinct vertices $t_1,t_2,t_3\in S\setminus R$.

If $s\in L$, then $d_{H_s}(s)=2$ and $(s,2)\notin E(H_s)$. Put $B=\{2,t_1,t_2,t_3\}$. None of the edges $(s,b)$, $b\in B$, belongs to $H_s$, and
\[
|B|=4=5-d_{H_s}(s)+1.
\]
Therefore, Lemma~\ref{lem:mobility} implies that some vertex in $B$ receives the leaf $9$. That vertex cannot be $2$, since this would produce $TS_5$. Hence one of $t_1,t_2,t_3$ belongs to $R$, which is a contradiction.

If $s\in\{3,5,7\}$, the same argument applies with $B=\{1,t_1,t_2,t_3\}$. The vertex $1$ cannot receive the leaf: it already has degree $5$ in $H_s$, so receiving $9$ would violate Lemma~\ref{lem:degree-bound}. Thus one of $t_1,t_2,t_3$ must receive the leaf, again contradicting their choice. It follows that $|R|\ge5$.
\end{proof}

We now show that if the edges \((2,c)\) and \((2,9)\),
for some \(c\in L\), belong to the same component,
then the partition $(\Gamma_1,\ldots,\Gamma_5)$ is weakly equivalent to one whose
fifth component is \(TS_5\).

\begin{lemma}\label{lem:d5-matching54}
Suppose that $\Gamma_5=H_a$ for some $a\in S$. If there exist $c\in L$ and $i\in\{1,2,3,4\}$ such that
\[
(2,c), \ \ \text{and} \ \ (2,9)\in E(\Gamma_i),
\]
then the partition $(\Gamma_1,\ldots,\Gamma_5)$ is weakly equivalent to one whose fifth component is $TS_5$.
\end{lemma}

\begin{proof}
If $a=c$, the involution $(c,2,9)$ transfers the leaf $9$ from $c$ to $2$, by Lemma~\ref{lem:leaf-switch}, and we are done. Suppose that $a\ne c$, and assume by contradiction that no weakly equivalent partition contains $TS_5$.

Let $(c,9)\in E(\Gamma_j)$. Since $(2,c)$ and $(2,9)$ both belong to $\Gamma_i$, we have $j\ne i$, or else these three edges would form a cycle in $\Gamma_i$. We write $\{i,j,k,l\}=\{1,2,3,4\}$. For each $r\in\{i,j,k,l\}$, let $n_r$ be the number of vertices $u\in S\setminus\{c\}$ for which $(c,u)\in E(\Gamma_r)$. The six edges from $c$ to $S\setminus\{c\}$ lie outside $H_a$ and are distributed among the four components $\Gamma_1,\ldots,\Gamma_4$. Thus
\begin{equation}\label{eq:d5-count54}
n_i+n_j+n_k+n_l=6.
\end{equation}

We first prove that
\begin{equation}\label{eq:d5-firstbound54}
n_i+n_j\le1.
\end{equation}

We consider the possible destinations of the leaf $9$,
which is initially attached to the vertex $a$. By
Corollary~\ref{cor:dminusone}, the leaf $9$ has at
least four destinations. The vertices eligible to
receive it are those in $S\setminus\{a\}$, together with one additional vertex: $2$ if $a\in L$, or $1$ if $a\in\{3,5,7\}$. The additional vertex cannot receive $9$, because the result would be $TS_5$ in the first case and would violate the degree bound in the second.

We also show that $c$ cannot receive $9$.
Indeed, if the involution $(a,c,9)$ transfers
the leaf $9$ from $a$ to $c$, then the edges
$(2,c)$ and $(2,9)$ remain in $\Gamma_i$.
By Lemma~\ref{lem:leaf-switch}, the involution
$(c,2,9)$ would then transfer the leaf $9$
from $c$ to $2$, producing $TS_5$, which is again a
contradiction.

Consequently, at least four vertices in
$S\setminus\{a,c\}$ receive the leaf $9$.
Together with its original attachment vertex
$a$, this gives at least five vertices in
$S\setminus\{c\}$ at which the leaf $9$
can be attached. Moreover, each of these
attachments can be obtained without involving
$c$ in any involution. Hence the component
memberships of all edges incident with $c$
remain unchanged.

Let $u$ be any of these five vertices. If $u\ne a$,
we transfer the leaf $9$ from $a$ to $u$ by applying
the involution $(a,u,9)$. If $u=a$, we make no change. We denote the resulting partition by
$(\Delta_1,\ldots,\Delta_5)$. Then $\Delta_5=H_u$,
and since the transfer does not involve $c$,
the edges incident with $c$ retain their original
component memberships. In particular,
\[
(c,9)\in E(\Delta_j),\quad \text{and} \quad
(2,c),(2,9)\in E(\Delta_i).
\]

We show that $(c,u)$ belongs to neither
$\Delta_i$ nor $\Delta_j$.
If $(c,u)\in E(\Delta_j)$, then $(c,u)$ and
$(c,9)$ belong to the same component. By
Lemma~\ref{lem:leaf-switch}, the involution
$(u,c,9)$ transfers the leaf $9$ from $u$ to $c$.
The involution $(c,2,9)$ then transfers it from
$c$ to $2$, producing $TS_5$, which is a contradiction.

Suppose instead that $(c,u)\in E(\Delta_i)$.
Applying the involution $(u,c,9)$ to the partition $(\Delta_1,\ldots,\Delta_5)$, it either transfers
the leaf $9$ to $c$, leading to the same
contradiction, or it leaves $H_u$ unchanged and
exchanges the component memberships of $(c,u)$
and $(c,9)$. In the latter case, the edges
\[
(2,c),\quad(c,9),\quad \text{and} \quad (9,2)
\]
all belong to the same component, forming a
cycle which is again a contradiction.

Therefore, for each of these five vertices $u$,
the edge $(c,u)$ belongs to neither $\Gamma_i$
nor $\Gamma_j$. Since $S\setminus\{c\}$ contains
exactly six vertices, at most one of the six
edges joining $c$ to these vertices can belong
to $\Gamma_i$ or $\Gamma_j$. Hence
\[
n_i+n_j\le1,
\]
as required.

By \eqref{eq:d5-count54} and \eqref{eq:d5-firstbound54},
we have $n_k+n_l\ge5$. After interchanging $k$ and $l$
if necessary, we may assume that $n_k\ge3$.

Recall that there are five vertices in $S\setminus\{c\}$
to which the leaf $9$ can be attached. Since $n_k\ge3$
and $S\setminus\{c\}$ contains six vertices, at least
two of these five vertices are joined to $c$ by an edge
of $\Gamma_k$. Choose one such vertex $r$.

If $r\ne a$, we apply the involution $(a,r,9)$ to
transfer the leaf $9$ from $a$ to $r$; otherwise,
we make no change. We denote the resulting partition
by $(\Delta_1,\ldots,\Delta_5)$. Then
$\Delta_5=H_r$, and since the transfer does not
involve $c$, we still have
\[
(c,r)\in E(\Delta_k),\qquad
(c,9)\in E(\Delta_j),
\]
and
\[
(2,c),(2,9)\in E(\Delta_i).
\]

We now apply the involution $(r,c,9)$ to the partition $(\Delta_1,\ldots,\Delta_5)$. If it transfers
the leaf $9$ from $r$ to $c$, then the involution
$(c,2,9)$ transfers $9$ from $c$ to $2$, producing
$TS_5$, which is a contradiction.

Otherwise, let $(\Theta_1,\ldots,\Theta_5)$ denote
the resulting partition. By Lemma~\ref{lem:leaf-switch},
the involution leaves $H_r$ unchanged and exchanges
the component memberships of $(c,r)$ and $(c,9)$.
Consequently,
\[
\Theta_5=H_r,\qquad
(c,9)\in E(\Theta_k),
\]
while
\[
(2,c) \ \ \text{and} \ \ (2,9)\in E(\Theta_i).
\]

Among the six edges $(c,u)$ with
$u\in S\setminus\{c\}$, only $(c,r)$ changes
component membership, moving from $\Gamma_k$
to $\Gamma_j$. Thus, if $n'_k$ denotes the number
of these edges belonging to $\Theta_k$, then
\[
n'_k=n_k-1\ge2.
\]

On the other hand, the argument used to establish
\eqref{eq:d5-firstbound54} applies to
$(\Theta_1,\ldots,\Theta_5)$, with $\Theta_k$
playing the role previously played by $\Gamma_j$.
Indeed, $(c,9)$ now belongs to $\Theta_k$, while
$(2,c)$ and $(2,9)$ still belong to $\Theta_i$.
The same argument therefore gives $n'_k\le1$,
which is a contradiction. This completes the proof.
\end{proof}

We are now ready to prove that the partition $(\Gamma_1, \cdots ,\Gamma_5)$ is weakly equivalent to a partition whose fifth component is $TS_5$.

\begin{lemma}\label{lem:d5-repeated54}
Suppose that $\Gamma_5=H_a$ for some $a\in L$. Then the partition $(\Gamma_1, \cdots ,\Gamma_5)$ is weakly equivalent to one whose fifth component is $TS_5$.
\end{lemma}

\begin{proof}

The five edges
\begin{equation}\label{eq2}
(2,4),\quad (2,6),\quad (2,8),\quad (2,10),
\quad \text{and}\quad (2,9)
\end{equation}
are precisely the edges incident with $2$ that do not
belong to $\Gamma_5$. Since $d_{\Gamma_5}(2)=4$,
each of the other four components contains at least
one of these five edges. Consequently, exactly one
of $\Gamma_1,\ldots,\Gamma_4$ contains two of them,
while each of the remaining three contains one.

If $(2,9)$ belongs to the component containing two
edges, then there exist $c\in L$ and
$i\in\{1,2,3,4\}$ such that
\[
(2,c) \ \ \text{and} \ \ (2,9)\in E(\Gamma_i).
\]
The conclusion then follows from
Lemma~\ref{lem:d5-matching54}.

We may therefore assume that there exist distinct
$p,q\in L$ and distinct $i,j\in\{1,2,3,4\}$
such that
\[
(2,p),(2,q)\in E(\Gamma_i),
\quad \text{and} \quad (2,9)\in E(\Gamma_j).
\]
Suppose, for a contradiction, that no partition
weakly equivalent to $\mathcal P$ contains a
component isomorphic to $TS_5$. Let $R$ be the
set defined in Lemma~\ref{lem:d5-reachable}.

We first show that neither $p$ nor $q$ belongs to
$R$. Suppose that $p\in R$. By
Lemma~\ref{lem:d5-reachable}, we can transfer the
leaf $9$ from $a$ to $p$ through a sequence of
involutions that keeps the eight edges of $C$
in the fifth component. We denote the resulting
partition by $(\Delta_1,\ldots,\Delta_5)$.
Since none of these involutions involves $2$,
we have
\[
\Delta_5=H_p,\qquad
(2,p),(2,q)\in E(\Delta_i),
\quad \text{and} \quad (2,9)\in E(\Delta_j).
\]

We now apply the involution $(p,2,9)$.
If it transfers the leaf $9$ from $p$ to $2$,
the resulting fifth component is $TS_5$,
contrary to our assumption. Otherwise, by
Lemma~\ref{lem:leaf-switch}, the involution
leaves $H_p$ unchanged and exchanges the
component memberships of $(2,p)$ and $(2,9)$.
Thus, in the resulting partition
$(\Theta_1,\ldots,\Theta_5)$, we have
\[
\Theta_5=H_p,\quad \text{and} \quad
(2,q),(2,9)\in E(\Theta_i).
\]
Lemma~\ref{lem:d5-matching54} then applies and gives a
contradiction. The same argument applies if
$q\in R$. Hence
\[
R\cap\{p,q\}=\varnothing.
\]

Since $R\subseteq S$, with $|R|\ge5$, and $|S|=7$,
we conclude that
\[
R=S\setminus\{p,q\}.
\]

We write $L\setminus\{p,q\}=\{r,s\}$, and let
\[
(2,r)\in E(\Gamma_k),\quad \text{and} \quad
(2,s)\in E(\Gamma_l).
\]
By the distribution of the five edges in
\eqref{eq2}, the indices $i,j,k,l$ are
pairwise distinct. Moreover,
$(p,q)\notin E(\Gamma_i)$, since otherwise
the edges $(2,p)$, $(2,q)$, and $(p,q)$
would form a cycle in $\Gamma_i$.
We now exclude the three remaining possibilities
for the component containing $(p,q)$.

Suppose first that $(p,q)\in E(\Gamma_j)$.
We apply the involution $(2,p,q)$ and denote
the resulting partition by
$(\Delta_1,\ldots,\Delta_5)$. Since
$(2,p)$ and $(2,q)$ belong to $\Gamma_i$,
while $(p,q)$ belongs to $\Gamma_j$,
the involution replaces one of $(2,p)$
or $(2,q)$ in $\Gamma_i$ by $(p,q)$.
Thus one of these two edges moves into
$\Delta_j$, which still contains $(2,9)$.
The fifth component remains $H_a$.
Lemma~\ref{lem:d5-matching54} therefore
applies to the resulting partition, giving
a contradiction.

Suppose next that $(p,q)\in E(\Gamma_k)$.
Since $r\in R$, we can transfer the leaf $9$
from $a$ to $r$, keeping the edges of $C$
in the fifth component. Denote the resulting
partition by $(\Delta_1,\ldots,\Delta_5)$.
Every intermediate attachment vertex belongs
to $R$, which contains neither $p$ nor $q$.
Furthermore, none of these transfers involves
$2$. Consequently,
\[
\Delta_5=H_r,\qquad
(p,q),(2,r)\in E(\Delta_k),
\]
while
\[
(2,p),(2,q)\in E(\Delta_i),
\quad \text{and} \quad (2,9)\in E(\Delta_j).
\]

We apply the involution $(2,p,q)$ to this
partition. Since $(2,p)$ and $(2,q)$
belong to $\Delta_i$, while $(p,q)$
belongs to $\Delta_k$, the involution
moves one of $(2,p)$ or $(2,q)$ into
the $k$-th component, leaving $H_r$
unchanged. After interchanging $p$ and $q$
if necessary, we may assume that
$(2,p)$ moves into that component.
Thus, writing $(\Theta_1,\ldots,\Theta_5)$
for the resulting partition, we have
\[
\Theta_5=H_r,\qquad
(2,p),(2,r)\in E(\Theta_k),
\quad \text{and} \quad (2,9)\in E(\Theta_j).
\]

We now apply the involution $(r,2,9)$.
If it transfers the leaf $9$ from $r$
to $2$, the resulting fifth component is
$TS_5$, which is a contradiction. Otherwise,
Lemma~\ref{lem:leaf-switch} implies that
the involution leaves $H_r$ unchanged
and exchanges the component memberships
of $(2,r)$ and $(2,9)$. Hence, in the
resulting partition
$(\Lambda_1,\ldots,\Lambda_5)$, we have
\[
\Lambda_5=H_r,\quad \text{and} \quad
(2,p),(2,9)\in E(\Lambda_k).
\]
Applying Lemma~\ref{lem:d5-matching54}
again gives a contradiction.

Finally, if $(p,q)\in E(\Gamma_l)$,
we repeat the preceding argument with
$s$ in place of $r$ and $\Gamma_l$
in place of $\Gamma_k$. This also
leads to a contradiction.

Therefore $(p,q)$ belongs to none of the five components \(\Gamma_1,\ldots,\Gamma_5\), which is a contradiction. Hence the partition $(\Gamma_1, \cdots ,\Gamma_5)$ is weakly equivalent to one whose fifth component is $TS_5$.
\end{proof}

Therefore, we have the following Corollary.

\begin{corollary}\label{lem:d5-no54-seven}
If the extremal component \(T\) satisfies \(\ell(T)=7\), then \(T\) cannot have branch-degree multiset $\{5,4\}$.
\end{corollary}

\begin{proof}
By Lemma~\ref{lem:d5-normal54}, after weak equivalence, we may assume that $T=H_a$ for some $a\in L$. Lemma~\ref{lem:d5-repeated54} gives a weakly equivalent partition whose fifth component is $TS_5$. Since $TS_5$ has eight leaves, this contradicts the maximality in \textup{(E1)}.
\end{proof}

Having excluded all three possible branch-degree multisets for an
extremal component with seven leaves, we can now establish our main result.

\begin{theorem}\label{thm:main}
Conjecture~\ref{conj:TSd} holds for $d=5$.
\end{theorem}

\begin{proof}
Let $\mathcal P$ be a cycle-free $5$-partition of $K_{10}$. Let $T$ be an extremal
component of $\mathcal P$ as in Definition~\ref{def:extremal-component}. Thus $\mathcal P$ is weakly
equivalent to a partition $(\Gamma_1,\ldots,\Gamma_5)$
with $\Gamma_5=T$.

By Theorem~\ref{thm:general-lower-bound} and
Lemma~\ref{lem:twinstar-characterization}, we have
\[
6\le\ell(T)\le8.
\]
Proposition~\ref{prop:d5-seven} and
Lemmas~\ref{lem:d5-no533-seven}, \ref{lem:d5-no443-seven} together with Corollary
\ref{lem:d5-no54-seven} exclude the cases
$\ell(T)=6$ and $\ell(T)=7$. Hence $\ell(T)=8$,
and Lemma~\ref{lem:twinstar-characterization}
implies that $T$ is isomorphic to $TS_5$.

Therefore, $\mathcal P$ is weakly equivalent
to a cycle-free $5$-partition containing a
component isomorphic to $TS_5$, proving
Conjecture~\ref{conj:TSd} for $d=5$.
\end{proof}
As a consequence of Theorem~\ref{thm:main}, we obtain
the following uniqueness result for the determinant-like
map $\det^{S^2}$ in dimension $5$.

\begin{corollary}\label{cor:det-uniqueness}
The determinant-like map
\(
\mathrm{det^{S^2}}:V_5^{\otimes45}\longrightarrow k
\)
is unique up to multiplication by a scalar.
\end{corollary}

\begin{proof}
By Theorem~\ref{thm:main} and the previously
established cases $d\le4$, the twin-star hypothesis
holds for every $1\le d\le5$. Hence
Conjecture~\ref{conj:Cd} holds for $d=5$ by
\cite{twin}. The uniqueness of $\det^{S^2}$
then follows from the equivalence established
in \cite{edge}.
\end{proof}

\section{Conclusion}\label{sec:conclusion}

One of the main difficulties in studying cycle-free
\(d\)-partitions of \(K_{2d}\) is that triangle involutions
affect the distribution of edges among different components,
making it difficult to control the entire partition
simultaneously. Our approach addresses this difficulty
by focusing on a single extremal component and examining
the remaining components only when necessary. This approach could potentially be extended to higher
dimensions, although, in its present form, it would require
treating each dimension separately. This leads to a natural
question:

\textbf{Question:} Can the twin-star hypothesis be established
in all dimensions through a uniform extremal approach?

One possible direction is to refine the extremal conditions
\textup{(E1)} and \textup{(E2)} of
Definition~\ref{def:extremal-component}, for instance,
by also maximizing the number of pendant paths attached
to a single branch vertex or by imposing additional
restrictions on the degrees of branch vertices.
Although such refinements were not necessary for $d=5$,
they may provide further structural information
toward a proof of the twin-star hypothesis in general.

\vspace{1cm}

\noindent\textbf{Acknowledgments.} The authors would like to thank Mihai D. Staic for insightful
discussions and valuable suggestions.

\end{document}